\documentclass[final,5p,times,twocolumn]{elsarticle}

\usepackage{amssymb}

\usepackage[ruled,linesnumbered]{algorithm2e}  
\usepackage{amsmath}  
\usepackage{subfigure}
\usepackage{booktabs}
\usepackage{rotating}
\journal{}
\biboptions{compress}
\graphicspath{{figure/}}
\begin{document}

\begin{frontmatter}



\title{A novel dual-stage algorithm for capacitated arc routing problems with time-dependent service costs}


\author[1]{Qingya Li}
\ead{liqy2020@mail.sustech.edu.cn}

\author[1]{Shengcai Liu}
\ead{liusccc@gmail.com}

\author[3,4]{Juan Zou}
\ead{zoujuan@xtu.edu.cn}
\author[1]{Ke Tang}
\ead{tangk3@sustech.edu.cn}

\address[1]{Guangdong Key Laboratory of Brain-Inspired Intelligent Computation, Department of Computer Science and Engineering, Southern University of Science and Technology, Shenzhen 518055, China}
\address[3]{The Key Laboratory of Hunan Province for Internet of Things and Information Security, Xiangtan University, Xiangtan 411105, China}
\address[4]{The Key Laboratory of Intelligent Computing and Information Processing, Ministry of Education, Xiangtan University, Xiangtan 411105, China}


\begin{abstract}
This paper focuses on solving the capacitated arc routing problem with time-dependent service costs (CARPTDSC), which is motivated by winter gritting applications.
In the current literature, exact algorithms designed for CARPTDSC can only handle small-scale instances, while heuristic algorithms fail to obtain high-quality  solutions.
To overcome these limitations, we propose a novel dual-stage algorithm, called MAENS-GN, that consists of a routing stage and a vehicle departure time optimization stage.
The former obtains the routing plan, while the the latter determines the vehicle departure time.
Importantly, existing literature often ignores the characteristic information contained in the relationship between the route cost and the vehicle departure time.
The most significant innovation in this paper lies in the exploitation of this characteristic information during the vehicle departure time optimization stage.
Specifically, we conduct a detailed analysis of this relationship under various scenarios and employ tailored methods to obtain the (approximately) optimal vehicle departure time.
Furthermore, we propose an improved initialization strategy that considers time-dependent characteristics to achieve better solution quality.
In addition to the modified benchmark test sets, we also experiment on a real-world test set.
Experimental results demonstrate that MAENS-GN can obtain high-quality solutions on both small-scale and larger-scale instances of CARPTDSC.
\end{abstract}



\begin{keyword}


Capacitated arc routing problem \sep
Time-dependent  service costs\sep
Vehicle departure time \sep
Memetic algorithm \sep
Combinatorial optimization

\end{keyword}

\end{frontmatter}


\section{Introduction}
\label{into}
The arc routing problem (ARP) is a complex planning and scheduling problem commonly found in operations research, transportation, and industrial applications \cite{jin2020planning}.
Among ARPs, the capacitated arc routing problem (CARP) is the most representative form with capacity constraints \cite{tang2009memetic}.
CARP is widely applied in numerous practical scenarios, such as winter gritting applications \cite{handa2006robust}, urban waste collection \cite{bodin1989design,lacomme2005evolutionary}, street sweeping \cite{bodin1978computer}, and electric meter reading \cite{stern1979routing}.
As a combinatorial optimization problem, CARP requires determining the least-cost routing plan for vehicles to serve all tasks subject to certain constraints \cite{tang2009memetic}.
CARP is generally described by a directed graph in which each task is associated with an arc.
Arcs with tasks that must be served are terms required arcs, and each required arc is associated with a service cost and a demand.
Arcs without tasks are called deadhead or travel arcs, which do not require service.
Each travel arc is always associated with a travel cost.
If certain information, such as service cost, travel cost, or demand, varies with time, then the problem becomes the time-dependent capacitated arc routing problem (TDCARP).
Among TDCARPs, the CARP with time-dependent service costs (CARPTDSC) is common in real-life scenarios. \\
\indent CARPTDSC is motivated by winter gritting or salting applications \cite{tagmouti2007arc, tagmouti2010variable}.
In this scenario, a fleet of fully loaded trucks, carrying deicing material (like salt), is dispatched to remove snow from critical roads in the urban road network.
The objective is to determine the routing plan and the vehicle departure time that minimizes the total cost while serving all tasks.
Each task is specified for a road from which a truck must remove snow by salting.
The service cost of each task is defined as the consumption of time or deicing material when serving the task.
However, this consumption always varies depending on the time at which the service begins for a particular road, due to factors like traffic, temperature, and weather conditions. 
Therefore, the time of beginning of service for a task affects its service cost, making the service cost time-dependent.
In \cite{tagmouti2007arc,tagmouti2011dynamic}, the relationship between the service cost of each task and its time of beginning of service  is modeled as a three-segment linear function, as depicted in Fig. \ref{fig0}(a).
In this case,  if the task is served at the optimal time point or interval (i.e., [\(bt\), \(et\)] in Fig. \ref{fig0}(a)), the service cost is minimal.
If service begins earlier or later than this time point or interval, the service cost of the task will increase dramatically.
For convenience, in \cite{tagmouti2010variable}, the relationship is modeled as a two-segment linear function, as shown in Fig. \ref{fig0}(b), which is a degenerate form of the three-segment linear function.
In this case, the service cost is lowest when the service starts at time 0, although this rarely happens in real life \cite{tagmouti2010variable}.
Based on the types of time-dependent functions, CARPTDSC can be classified into two-segment linear function problems (abbreviated as 2LP) and three-segment linear function problems (abbreviated as 3LP).

For the 2LP, the variable neighborhood descent heuristic (VND) \cite{tagmouti2010variable} has been proposed, but it has been observed to be less effective.
On the other hand, for the 3LP, an exact algorithm~\cite{tagmouti2007arc} has been proposed, which however is only capable of handling small-scale problem instances (number of tasks \(\leq 40\)).
Motivated by these limitations, this paper investigates two research questions: 
\begin{description}
	\item[RQ1] For the 2LP, the quality of the approximate solutions obtained by the existing algorithm is not high, how to further enhance the quality of the solutions?
	\item[RQ2] For the 3LP, the existing algorithm is limited to solving small-scale instances, how to effectively address larger-scale instances?
\end{description}

To address the aforementioned research questions, we propose a novel dual-stage algorithm.
Given that  CARPTDSC involves  two types of decision variables, namely the routing plan $X$ and the vehicle departure time $ T $,  the two stages of the proposed algorithm align with solving for these decision variables, respectively. Consequently, the dual stages  are denoted as the routing stage and the vehicle departure time optimization stage, respectively. In the routing stage, the primary focus is on obtaining a routing plan \(X\) by modifying the memetic algorithm with extended neighborhood search (MAENS) \cite{tang2009memetic}. On the other hand, the vehicle departure time optimization stage concentrates  on determining  the vehicle departure time corresponding to each route in \(X\) using the golden section search (GSS) \cite{OVERHOLT1967FIBONACCI} or negatively correlated search (NCS) \cite{tang2016negatively}. As a result, the proposed algorithm is termed MAENS-GN. The contributions made in this paper are briefly summarized as follows.
   \begin{itemize}
   	
   	\item When optimizing for vehicle departure time, existing works have ignored the characteristic information embedded  in the relationship between the route cost and the vehicle departure time. To exploit this information to improve performance,  we conduct a detailed analysis of the relationship  across different scenarios. Based on the relationships analyzed in different scenarios, tailored approaches are adopted.
    \item  An improved initialization strategy that considers time-dependent characteristics is proposed to enhance the effectiveness.
      	
   	\item Experiments were conducted on both modified standard test sets and a test set derived from real-world applications in Jingdong Logistics. The experimental results demonstrate that  MAENS-GN yields positive effectiveness, indicating that this algorithm has addressed two research questions well.
   \end{itemize}

The rest of this paper is organized as follows. Section \ref{backrw} introduces the background and related work. Section \ref{algo} describes the proposed algorithm, i.e., MAENS-GN. After that, experimental studies are presented in Section \ref{expe}. Finally, Section \ref{concf} presents the conclusion and future work.


\section{Background and related work}
\label{backrw}
In this section, the definition of CARPTDSC is first introduced. Then, the static algorithm MAENS \cite{tang2009memetic}, especially its initialization strategy, is described. Finally, we introduce the related work.
\subsection{Problem definition}
\label{defi}
$G = (V, A) $ is a directed graph where $ V $ is the set of vertices and $ A $ denotes the set of arcs. Each arc $e \in A$, which can be represented as \( e=\langle i, j \rangle \), and \(i, j \in V\). \(i\), \(j\) represent the tail node and head node of arc \(e\), respectively.  $ A $ can be partitioned  into two subsets, namely $ A_1 $ and $ A_2 $. $ A_1 \subseteq A $ denotes the set of required arcs, also known as the set of tasks, each of which must be served. Each required arc is assigned a unique ID, which is set to the positive integer.
$ A_2 \subseteq A $ represents the set of travel or deadhead arcs, each of which does not need to be served and is used only for a traveling action. 
Each required arc $ re \in A_1 $  is associated with a demand $ d(re) $, a length $ len_{re} $, a deadhead or travel time $ dt_{re} $, a service time $ st_{re} $, a deadhead or travel cost $ DC_{re} $, and a time-dependent piecewise linear service cost function $SC(re, T_{re})  $, where $T_{re}  $ denotes the time of beginning of service on required arc $ re $. Each travel arc $de \in A_2  $ is associated with a length $ len_{de} $, a travel time $ dt_{de} $, and a travel cost $ DC_{de} $.  

\begin{figure}[htbp]
	\centering
	   \includegraphics[width=0.5\textwidth]{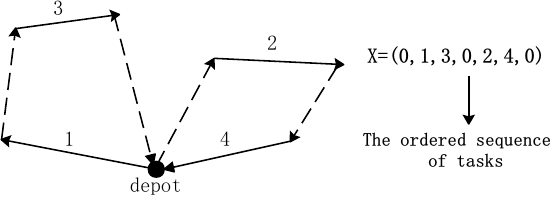}
	\centering
	\caption{Illustration of the routing plan representation. Each solid line and each dotted line represent a required arc (task) and a travel arc, respectively.}
	\label{fig_r}
\end{figure}

$ K $ identical vehicles with capacity $ Q $ are available to serve these required arcs. Initially, all vehicles are stationed at the central depot from which each vehicle starts a single route at a certain vehicle departure time. The route starts and ends in the depot. For the sake of convenience, the depot is used as a dummy required arc (task), represented by the unique ID 0. The value of each attribute of the depot is set to 0, such as service cost, service time and demand, etc.
Notably, since no waiting time is permitted  on the route, the arrival time of the task is equal to its time of beginning of service.
The objective is to get a feasible solution with minimum total cost.  In \cite{tagmouti2007arc} \cite{tagmouti2010variable}, it is noted that in  winter gritting applications, if the vehicle departure time is different, the cost in time will be different. Consequently,  in this context, the total cost incurred  is specified as the total time expended, which means that the service cost and the travel cost are equal to the service time and the travel time, respectively. Therefore, in this paper, the time-dependent service cost is actually the time-dependent service time.
The mathematical form of the objective function is as follows:
\begin{equation}\label{equ:1}
min \quad C(X, T)=\sum_{k=1}^{m} C(R_k, t_k),
\end{equation}
where have two types of decision variables, denoted as $ X $ and $ T $,   representing the routing plan and the vehicle departure time, respectively. Therefore, the feasible solution \(S\) consists of two parts, i.e., \(X\) and \(T\), which can be expressed as \(S=(X, T)\). 
Specifically, $ X=(R_1, R_2, ..., R_m) $ and $ T=(t_1, t_2, ..., t_m) $,  where $ m $ indicates  the number of routes or vehicles. Each $ R_k $ represents the route served and traveled by vehicle $k$, expressed as $  R_k=(re_0^k, re_1^k, re_2^k,..., re_{l_k}^k,$ $ re_{l_k+1}^k)$. In essence, $ R_k $ comprises  depot 0, required arcs $(re_1^k, re_2^k,..., re_{l_k}^k)  $ and travel arcs $(de_1^k,$ $ de_2^k, ..., de_{p_k}^k)$. 
These travel arcs result from the shortest paths between every two neighboring (dummy) required arcs in the route. The shortest paths are gotten  through Dijkstra's algorithm \cite{dijkstra1959note}. Fig. \ref{fig_r} presents a simple illustration of the routing plan representation. The routing plan \(X=(0,1,3,0,2,4,0)\). In \(X\),  0 is the separator, denoting the depot. And each \(r \in (1, 3, 2, 4)\) represents a task. Notably, every sequence of tasks enclosed by two separators (including separator 0) represents a route. So \(X\) comprises two routes, i.e., \((0,1,3,0)\) and \((0,2,4,0)\). \\
\indent Additionally, $C(R_k, t_k)  $ denotes the cost of the route $ R_k $ at vehicle departure time \(t_k\), with a mathematical formula as follows:
\begin{equation}\label{equ:2}
C(R_k, t_k)=\sum_{i=1}^{l_k} SC(re_i^k, T_{re_i^k}(t_k))+\sum_{i=1}^{p_k} DC_{de_i^k},
\end{equation}
where $ l_k $ and $ p_k $ represent the numbers of the required arcs (except the dummy required arc, i.e., depot) and the  travel arcs in the route $ R_k $, respectively, and $ k \in \{1, 2, ..., m\} $. $ SC(re_i^k, T_{re_i^k}(t_k)) $ denotes the service cost of the required arc $ re_i^k $ with its the time of beginning of service $ T_{re_i^k}(t_k) $, where $ i \in \{1, 2, ..., l_k\} $.  $ T_{re_i^k}(t_k) $ can be further expressed as $T_{re_i^k}(t_k)=t_k+\sum_{h=0}^{i-1} st_{re_h}^k +\sum_{h=0}^{i-1} dt_{sp_h}^k $, where $st_{re_h}^k$   denotes the service time of the $ h $-th required arc and $ dt_{sp_h}^k $ is the travel time of the shortest path $ sp_h $ between the $ h $-th and $ (h+1)$-th required arcs and $ h \in \{0, 1, ...,i-1\} $.  Additionally, $ DC_{de_i^k} $ represents  the travel cost of the travel arc $ de_i^k $, where $ i \in \{1, 2, ...,p_k \} $.

In addition to the objective function, CARPTDSC has some constraints whose mathematical formulas are as follows.
\begin{equation}\label{equ:3}
re_0^k=re_{l_k+1}^k=0,
\end{equation}
\begin{equation}\label{equ:4}
re_i^k \neq re_{i'}^{k'}, \forall (k, i) \neq (k',i'),
\end{equation}
\begin{equation}\label{equ:5}
re_i^k  \neq inv(re_{i'}^{k'}), \forall(k,i)\neq (k',i'),
\end{equation}
\begin{equation}\label{equ:6}
\sum_{k=1}^{m} l_k=|A_1|,
\end{equation}
\begin{equation}\label{equ:7}
\sum_{i=1}^{l_k} d(re_i^k)<Q,
\end{equation}
\begin{equation}\label{equ:8}
0 \leq T_{re_i^k}(t_k) \leq T.
\end{equation}

In the Eqs. (\ref{equ:3})-(\ref{equ:8}), \(k \in \{1, 2, \ldots, m\}\). In the Eqs. (\ref{equ:4})-(\ref{equ:5}), \(i \in \{ 1, 2, \ldots, l_k\} \), and \(i \in \{0, 1, \ldots, l_k, l_k+1 \} \) in the Eq. (\ref{equ:8}). 
Eq. (\ref{equ:3}) represents the constraint that each route must begin and end at the depot.  Eqs. (\ref{equ:4})-(\ref{equ:6}) denote the constraint that each required arc can be served only once. In Eq. (\ref{equ:5}), $ inv(re_{i'}^{k'}) $ represents the arc $re_{i'}^{k'}  $ with the reverse direction of itself (if $ inv(re_{i'}^{k'}) $ exists). $|A_1|$ in Eq. (\ref{equ:6}) denotes the number of all required arcs in the set  of required arcs $A_1$. Eq. (\ref{equ:7}) is the capacity constraint which means that the total demand of each route doesn't exceed the capacity of the vehicle.  The last constraint is the time constraint, as shown in Eq. (\ref{equ:8}), which indicates that the time of beginning of service on each (dummy) required arc  must be within the planning horizon $[0, T]$.

\subsection{MAENS}
\label{maens}
MAENS \cite{tang2009memetic} is a state-of-the-art method for solving static CARP. It includes initialization, crossover, three local search operators, and a novel large-neighborhood operator. This section mainly details the initialization strategy of MAENS.
\subsubsection{Initialization strategy of MAENS}
\label{maensini}
The initialization strategy in MAENS is to generate the initial population with $psize$ distinct  individuals. This process of generating an individual is termed individual initialization strategy here, whose main purpose is to generate an initial routing plan that contains all tasks.


The individual initialization strategy is a constructive heuristic strategy, and its process of constructing the routing plan is the process of constructing routes one by one. When constructing a route,  the route is first initialized to empty, and then 0 (i.e., depot) is added to the route. Then one task is  added iteratively to the route. When selecting the task to be added, the task nearest to the end of the current route (i.e., the head node of the last task) is selected from all the unserved and feasible tasks. If multiple tasks meet the selection condition at the same time, the tie-breaking rule selected at this time is to randomly select one from multiple tasks that meet the condition. When there is no feasible task to select, 0 is added to the route. At this time, the process of constructing a route ends.

\subsection{Related work}
\label{rw}
In \cite{golden1981capacitated}, CARP was introduced  for the first time, and it was proved that CARP is NP-hard. Due to the limitations of exact algorithms in handling only small-scale instances, the research community has increasingly embraced heuristic and meta-heuristic approaches.  Prominent  heuristic algorithms include  construct-strike \cite{christofides1973optimum}, path-scanning \cite{golden1983computational}, and augment-merge \cite{golden1981capacitated}.
And many meta-heuristic methods are also applied, including tabu search \cite{hertz2000tabu},
variable neighborhood descent \cite{hertz2001variable}, guided local search \cite{beullens2003guided}, evolutionary algorithm \cite{handa2006robust,handa2006robustso}, memetic algorithm (or hybrid genetic algorithm) \cite{tang2009memetic,lacomme2004competitive,vidal2017node}.

Compared to ARP, vehicle routing problem (VRP) or node routing problem is a more extensively  studied routing problem,  commonly applied in logistics. Due to the extensive research on VRP, some ARPs are transformed into VRPs for resolution \cite{longo2006solving,pearn1987transforming}, including CARPTDSC \cite{tagmouti2007arc}. And among  various VRPs, TDVRP is the one that has the same kind of characteristics as the time-dependent arc routing problem (TDARP), i.e., time-dependent. Therefore, understanding the solving algorithms for TDVPR is also very helpful to address  TDARPs.  
In TDVRP, the travel time is generally time-dependent, influenced by factors such as the distance between two customers and the time of day. Numerous  TDVRPs have been formulated as mixed integer programming models \cite{malandraki1992time,chen2006real}. Algorithms for dealing with TDVRP include exact algorithms, heuristic algorithms, and meta-heuristic algorithms. Among them, the exact algorithms include branch-and-price method \cite{dabia2013branch} and dynamic programming \cite{kim2016solving,soysal2017simulation}. The heuristic algorithm has mathematical-programming-based heuristic \cite{malandraki1992time}. Additionally, meta-heuristic algorithms encompass genetic algorithm \cite{haghani2005dynamic}, and ant colony algorithm \cite{balseiro2011ant,yao2015vehicle}.

Although some ARPs can be transformed into  VRPs  for resolution, it is advantageous to design algorithms specifically for ARPs. In the original form of the problem, it is easier to exploit the characteristics of the problem \cite{tagmouti2011dynamic}. In addition, the transformation from VRP to ARP increases the dimension of the problem and reduces the attempt to use the VRP algorithms to solve ARP \cite{tong2021hybrid}. Regarding TDARPs, some algorithms have been designed.
Sun et al. \cite{sun2011new} proposed  a new integer programming formulation for the Chinese postman problem with time-dependent travel times, which significantly reduces  the size of the circuit formulation. Black et al. \cite{black2013time} studied a  time-dependent prize-collecting arc routing problem in which travel time changes with the time of day. Then, two metaheuristic algorithms were proposed to address this problem.  Nossack et al. \cite{nossack2017windy} proposed the windy rural postman problem with a time-dependent zigzag option, which was modeled as (mixed) integer programming formulation and then an exact algorithm was used to solve the problem. Calogiuri et al. \cite{calogiuri2019branch} mainly dealt with the time-dependent rural postman problem where the travel time of certain arcs varies in the planning horizon. Branch-and-bound algorithm was developed to solve this problem. 

Jin et al. \cite{jin2020planning} studied a capacitated arc routing problem with time-dependent penalty costs, which stems  from the planning of garbage collection services.  Penalty costs are related to parking patterns and service periods. The problem was formalized  as a mixed integer linear model and
a problem-specific intelligent heuristic search approach was proposed for its resolution. Ahabchane et al. \cite{ahabchane2020mixed} dealt with  the mixed capacitated general routing problem with time-dependent demands. CPLEX was used to solve small-scale instances and slack induction by string removals metaheuristic method was developed to address large-scale instances. Additionally, Vidal et al. \cite{vidal2021arc} studied a capacitated arc routing problem with time-dependent travel times and paths, marking the first exploration of CARP with time-dependent travel times at the network level.  Branch-and-price algorithm and hybrid genetic search were proposed for solving such problems.
\begin{algorithm}[!htb]  
	\caption{The general framework of MAENS-GN}  
	\label{alg:3}  
	\KwIn{A TDCARP instance $inst$\; } 
	\KwOut{The routing plan $X_{bf}$ and the vehicle departure time $Dt$\;}  
	Get the routing plan $X_{bf}  $ by MAENS($inst$)\;
        Get the instance type \(intp\) from \(inst\)\;
        \tcp{ Based on the analyzed relationship between the route cost and vehicle departure time under various scenarios (2LP, 3LP (slope \( |k| \leq 1\) and \(|k|>1\))), vehicle departure time begins to be optimized.}
       
        \If{\(intp == 2LP\)}
        {
         n $\leftarrow$ the number of routes in $X_{bf}$\;
         \For{$i\leftarrow 1:n$ }
         {
         $ Dt \leftarrow Dt \cup 0 $\;
         }
        }
        \ElseIf{\(intp == 3LP\) }
        {
        
	\If{$|k| \leq 1$}
	{
		$Dt \leftarrow \boldsymbol{GSS(X_{bf},\ 0,\ T)}$\;
	}
	\ElseIf{$|k| > 1$}
	{
		$Dt \leftarrow \boldsymbol{NCS(X_{bf},\ 0,\ T)}$\;
	}
        }
	\Return $X_{bf}$ and $Dt$
\end{algorithm}

\begin{figure}[htbp]
	\raggedright
	\noindent \makebox[\textwidth][l]{
	\subfigure[\textbf{ Three-segment linear function}]{
		\begin{minipage}[t]{0.5\linewidth}
			\raggedright
			\includegraphics[width=1.8in]{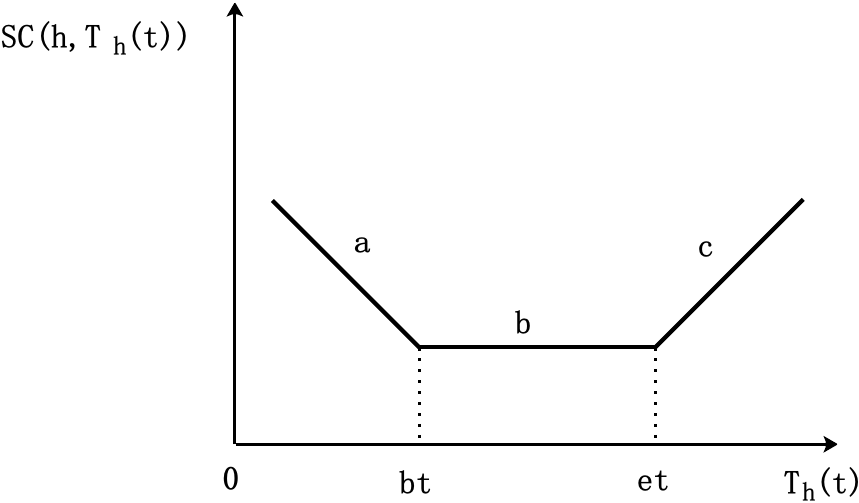}
			
		\end{minipage}%
	}%
	\subfigure[ \textbf{Two-segment linear function}]{
		\begin{minipage}[t]{0.5\linewidth}
			\centering
			\includegraphics[width=1.8in]{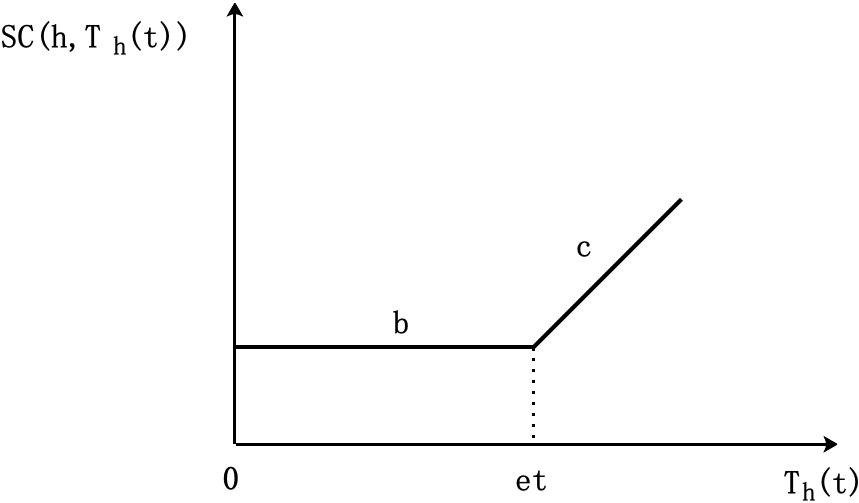}
			
		\end{minipage}%
	}%
	
}
\centering
	\caption{Schematic diagram of two types of  linear functions between the service cost of task $ h $ and its time of beginning of service. \(a\), \(b\), and \(c\)  represent the decreasing,  unchanged, and increasing segments, respectively.}
	\label{fig0}
\end{figure}
Despite several studies on various types of TDARPs, limited research has been dedicated solely to CARPTDSC. Tagmouti et al. \cite{tagmouti2007arc} first introduced time-dependent service costs to CARP which originates from winter gritting applications. 
Since the relationship between the service cost of each task and its time of beginning of service was modeled as a three-segment linear function, the type of CARPTDSC is 3LP. 
This problem was initially  transformed into a node routing problem and addressed through a column generation approach. However, the exact algorithm can only address the small-scale instances in which the number of tasks is only up to 40. To tackle  larger-scale instances, Tagmouti et al. \cite{tagmouti2010variable} proposed  a variable neighborhood descent heuristic (VND). For simplicity, a two-segment linear function was employed to represent  the relationship between the service cost of the task and its time of beginning  of service, which is a  degenerate form of three-segment linear function. Therefore, the type of CARPTDSC is 2LP, in which the route has the lowest cost if the corresponding vehicle starts service from time 0. In this case, the need for further departure time determination is eliminated and the problem's complexity is reduced.
What's more,  VND was less effective, which fails to reach optimal solutions in all test instances. To incorporate more practical considerations, Tagmouti et al. \cite{tagmouti2011dynamic} introduced a dynamic factor to CARPTDSC, i.e., changes in the optimal service time interval due to updates from weather forecasts. An adapted VND was used to address this problem. In the above papers \cite{tagmouti2007arc, tagmouti2010variable, tagmouti2011dynamic}, it is mentioned that the cost in time is time-dependent. However, the service time of the task is still assumed to be constant. Therefore, this is a bit unrealistic.

\section{The proposed dual-stage algorithm}
\label{algo}
The proposed dual-stage algorithm is termed MAENS-GN, which is divided into two stages, i.e., the routing stage and the vehicle departure time optimization stage. The general framework is described in Algorithm \ref{alg:3}.
In the routing stage,  MAENS \cite{tang2009memetic} was used to get the routing plan, i.e., $X_{bf}$ (line 1).  When using MAENS to get  $X_{bf}$, the evaluations need to be modified into time-dependent evaluations, which will not be described here. Furthermore, the initialization strategy of MAENS has been adapted  to novel initialization strategy. After  $X_{bf}$ is obtained, the vehicle departure time \(Dt\) begins to be optimized based on the analyzed relationship between the route cost and the vehicle departure time under various scenarios in the vehicle departure time optimization stage.
If the problem is 2LP,  time 0 is the optimal vehicle departure time (lines 3-8).  If the problem is 3LP,  suitable  methods are employed to search the (approximately) optimal departure time from the planning horizon [0, \(T\)] according to the absolute value (i.e., \(|k|\)) of the slope (lines 9-16). MAENS-GN introduces two key innovations: the optimization of vehicle departure time and the novel initialization strategy. The subsequent sections provide a detailed description of these two components.

\subsection{Optimization of the vehicle departure time}
\label{algosec:2}
 This optimization of the vehicle departure time is divided into two parts below, namely, the  analysis of the relationship between the route cost and the vehicle departure time, and the  methods for optimizing  the vehicle departure time.
\subsubsection{Detailed analysis of the relationship between the route cost and the vehicle departure time}
\label{math}
\begin{figure*}[htbp]
	\raggedright
	\noindent \makebox[\textwidth][l]{
	\subfigure[\textbf{ $h_1  $}]{
		\begin{minipage}[t]{0.3\linewidth}
			\raggedright
			\includegraphics[width=1.8in]{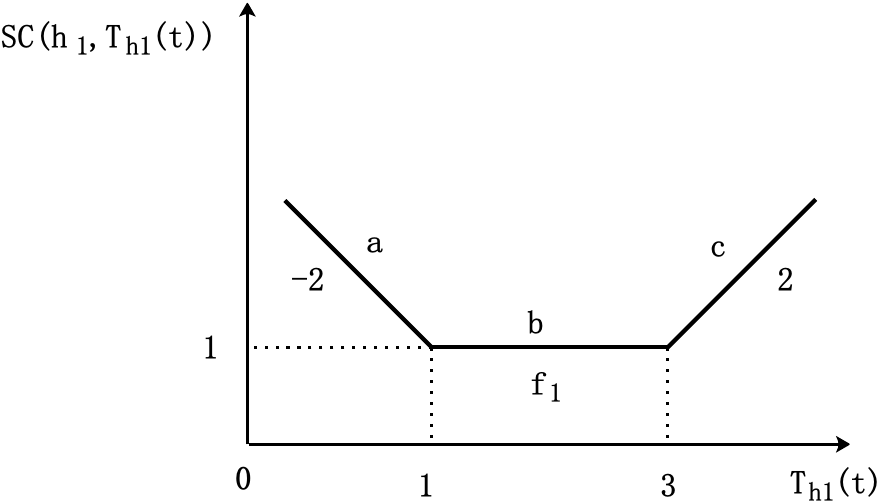}
		
		\end{minipage}%
	}%
	\subfigure[ \textbf{$ h_2 $}]{
		\begin{minipage}[t]{0.3\linewidth}
			\centering
			\includegraphics[width=1.8in]{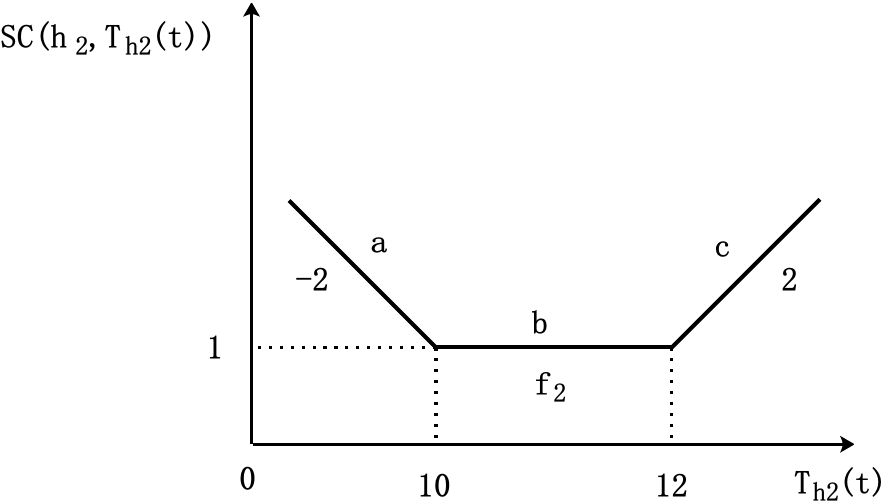}
		
		\end{minipage}%
	}%
	\subfigure[\textbf{$ h_3 $}]{
		\begin{minipage}[t]{0.3\linewidth}
			\centering
			\includegraphics[width=1.8in]{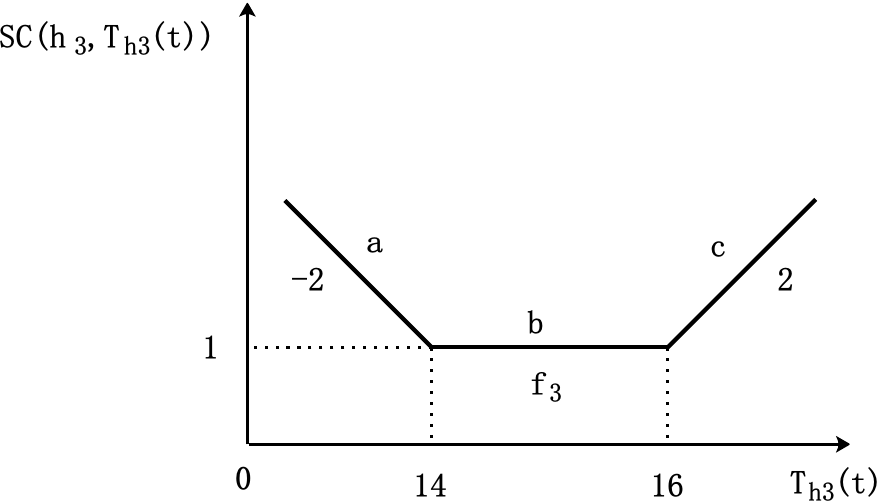}
			
		\end{minipage}
	}%
}

	\centering
	\caption{Three-segment linear functions of the  service cost of the three tasks ($ h_1 $, $ h_2 $ and $ h_3 $) in the route and their  time of beginning of service.}
	\label{fig4}
\end{figure*}

When analyzing the relationship between the route cost and the vehicle departure time, the travel cost is not calculated because it has no relationship with the vehicle departure time. It is worth noting that in this paper, the time-dependent service cost is specified for the time-dependent service time.
We performed a detailed analysis of the relationship between the route cost and the vehicle departure time for both the 2LP and 3LP (slope $|k|<1  $ and $|k|>1  $) scenarios. Since the relationship  under \(|k|=1\) in the 3LP is similar to that under the  \(|k|<1\), it will not be analyzed here. For convenience, in the 3LP, it is assumed that the absolute values of the slopes of  three-segment linear functions of all tasks in an instance are the same, as done in \cite{tagmouti2007arc}. For example, the absolute value of the slope is 2, which means the slopes of all tasks  are 2 and -2, as shown in Fig. \ref{fig4}.
\paragraph{Detailed analysis of the relationship between the route cost and the vehicle departure time in the 2LP}\hspace{0pt}

1. When there is only one task, i.e., $h$, in the route,  the relationship between the route cost and the vehicle departure time is consistent with the non-decreasing two-segment linear function of the task,  as illustrated  in Fig. \ref{fig0}(b). 

\begin{figure*}[htbp]

	\raggedright
	\noindent \makebox[\textwidth][l]{
	\subfigure[\textbf{ (similar) unimodal  }]{
		\begin{minipage}[t]{0.3\linewidth}
			\raggedright
			\includegraphics[width=1.8in]{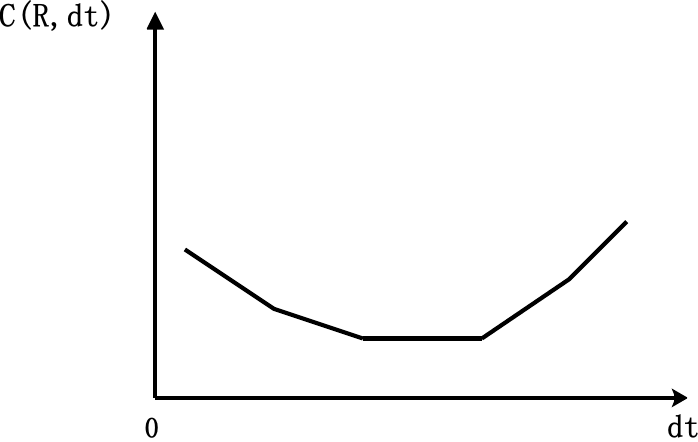}
			
		\end{minipage}%
	}%
	\subfigure[ \textbf{non-decreasing }]{
		\begin{minipage}[t]{0.3\linewidth}
			\centering
			\includegraphics[width=1.8in]{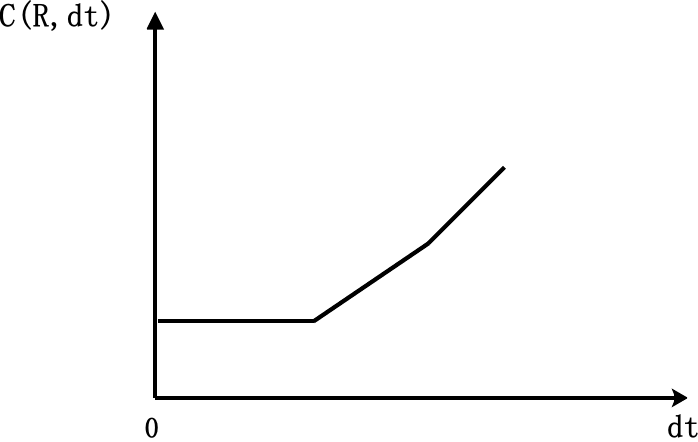}
			
		\end{minipage}%
	}%
	\subfigure[\textbf{non-unimodal }]{
		\begin{minipage}[t]{0.3\linewidth}
			\centering
			\includegraphics[width=1.8in]{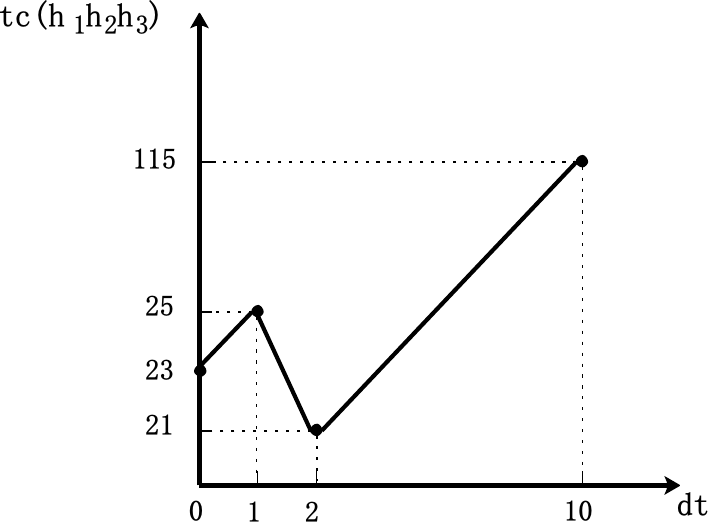}
			
		\end{minipage}
	}%
}

	\centering
	\caption{Schematic diagram of the relationship between the route cost \( C(R, dt)\) or \(tc(h_1h_2h_3)\) and vehicle departure time \(dt\).}
	\label{fig_m}
\end{figure*}

2. When there are multiple tasks in the route, the service cost of each task has a non-decreasing relationship with its  time of beginning of service. As the vehicle departure time increases, the time of beginning of service of each task will also increase. Consequently, the service cost of each task will also be non-increasing.  The route cost is equal to the sum of the service costs of all tasks in the route. Therefore, the route cost will have a non-decreasing relationship with the vehicle departure time, as shown in Fig. \ref{fig_m}(b).

Hence, whether there are one or more tasks in a route, the route cost has a non-decreasing relationship with the vehicle departure time. Consequently, time 0 is the optimal vehicle departure time for this route in the 2LP. This is also the reason why in \cite{tagmouti2010variable}, the algorithm does not need to solve for the vehicle departure time.
\paragraph{Detailed analysis of the relationship between route cost and vehicle departure time in the 3LP with the slope $|k|<1 $}\hspace{0pt}
\label{mathks1}
\indent When \(|k|<1\) in the 3LP, the analysis needs to be divided into two stages. Firstly, it is necessary to analyze the relationship between the time of beginning of service  of each task in the route and the vehicle departure time. Secondly, the relationship between the  route cost and the vehicle departure time is analyzed. In the first stage, it has been observed that as the vehicle departure time increases, the time of beginning of service of each task in the route also increases. The main analysis process is as follows. \\
\indent Assume that there are \(n\) tasks in the route. When the vehicle departure time \(dt\) increases by \(\Delta t\) (\(\Delta t\) is a  sufficiently small positive real number), if the increase of time of beginning of service of the \(i\)-th (\(1 \leq i \leq n\)) task in the route is to be minimal, the time of beginning of service  of the tasks before \(i\)-th task  in the route need to be in the decreasing segments (i.e., segment \(a\) in Fig. \ref{fig0}(a)). In this case, when \(dt\) increases by  \(\Delta t\), the time of beginning of service  of the first task increases by \(\Delta t\), and its service time (i.e., service cost) will decrease by \(|k|\Delta t\). Consequently, the time of beginning of service of the second task increases by \((1-|k|) \Delta t\), and its service time will decrease by \(|k|(1-|k|)\Delta t\). Therefore, the time of beginning of service  of the third task will increase by \((1-|k|)^{2} \Delta t\). It can be deduced that the time of beginning of service of the \(i\)-\(th\) task will increase by \((1-|k|)^{i-1} \Delta t\).  Because \(|k|<1\), \(0<(1-|k|)<1\).
Therefore,  among all the tasks in the route, the last task (i.e., the \(n\)-th task)  experiences  the least increase in time of beginning of service, which  increases by \((1-|k|)^{n-1} \Delta t\).  This value is still positive. 
Why is it that when the time of beginning of service of the tasks before \(n\)-th task in the route is all within the decreasing segments,
the increase of the time of beginning of service of the \(n\)-th task is minimal as \(dt\) increases by \(\Delta t\)? Below is an example.
Among the \(n\) tasks in the route,  the time of beginning of service of only the \(i\)-th (\(i<n\)) task is not within the decreasing segment, but within the increasing segment (i.e., segment \(c\) in Fig. \ref{fig0}(a)). At this time, the time of beginning of service of the \(n\)-th task is \((1+|k|)(1-|k|)^{n-2}\), which is greater than \((1-|k|)^{n-1}\). Therefore, as \(dt\) increases, the time of beginning of service of each task in the route also increases. The analysis of the second stage is as follows. \\
\indent When \(dt\) is 0, the time of beginning of service of the tasks in the route lies within different segments (e.g., segments \(a\), \(b\), or \(c\) in Fig. \ref{fig0}(a)) of the three-segment linear functions. Therefore,
when \(dt\) increases by \(\Delta t\) from 0 (i.e., \(dt=\Delta t\)), the service costs of the tasks in the route may be decreasing,  unchanged, or increasing.  Assume that  the change of the route cost at this time  is \(tc_0\). As \(dt\) continues to increase with the size of \(\Delta t\), the time of beginning of service  of all tasks in the route will also continue to increase, gradually approaching the increasing segment (i.e., segment \(c\) in Fig. \ref{fig0}(a)) or continuing to increase in the increasing segment. In the end, the time of beginning of service  of all tasks will fall within the increasing segment. Consequently, if \(tc_0<0\), as \(dt\) increases continually, the change in route cost will gradually become greater than 0. During this period, the route cost will initially  decrease (then remain unchanged), and finally increase, showing  a (similar) unimodal trend as shown in Fig. \ref{fig_m}(a). Conversely, if \(tc_0 \geq 0\), with the continued increase of \(dt\), the relationship between the route cost and the vehicle departure time will show a non-decreasing trend, as shown in Fig. \ref{fig_m}(b). \\
\indent Therefore, in the 3LP with slope \(|k|<1\), the relationship between the route cost and the vehicle departure time  exhibits  a (similar) unimodal or non-decreasing trend.

\paragraph{Detailed analysis of the relationship between route cost and vehicle departure time in the 3LP with slope $|k|>1 $}\hspace{0pt}
\label{mathkl1}

When  $ |k|>1 $, the relationship between the route cost and the vehicle departure time no longer meets the (similar) unimodal property. This is mainly because when  \(dt\) increases by $\Delta t$,  the time of beginning of service of the task in the route may not  increase, but appears to decrease. An illustrative example is provided  below.

A route contains three tasks (i.e., $ h_1 $, $ h_2 $, $ h_3 $), and their three-segment linear functions (i.e., $ f_1 $, $ f_2 $, $ f_3  $) are shown in Fig. \ref{fig4}, and the absolute value of the slopes is 2. The ranges of segments $  b$ of  $ f_1 $, $ f_2 $, and $ f_3 $ are [1, 3], [10, 12] and [14, 16], respectively. Among the three functions, the minimum  service costs of the tasks are 1. The route can be represented by the sequence of tasks $ h_1h_2h_3 $, with the depot positioned  at the tail node of $ h_1 $ and the head node of \(h_3\). The three tasks are directly connected, signifying that the head node of $ h_1 $ is equal to the tail node of $ h_2 $, and the head node of $ h_2 $ is equal to the tail node of $ h_3 $. Let  $T_h(dt)$   (where $h \in \{h_1, h_2, h_3\}$)  represents  the time to arrive at task $h$, and $SC(h, T_h(dt))$ represents the service cost of task $ h  $ at \(T_h(dt)\). Assume that the time $ \Delta t $ by which \(dt\) increases each time is 1.

When $ dt=0 $,  $ T_{h_1}(dt=0)=0 $. At this time, according to $ f_1 $, $SC(h_1, 0)=(1-0)*2+1=3 $. And $ T_{h_2}(0)=SC(h_1, 0)+T_{h_1}(0)=3 $. According to  $f_2  $, $SC(h_2, 3)=15  $. And $ T_{h_3}(0)=T_{h_2}(0)+SC(h_2, 3)=18 $. According to $f_3  $, $SC(h_3, 18)=5  $. Let $ tc(x) $ represent the sum of service costs of all tasks in the sequence of tasks $ x $. Therefore, $ tc(h_1h_2)=SC(h_1,0)+SC(h_2,3)=18 $, and $tc(h_1h_2h_3)=SC(h_1,0)+SC(h_2,3)+SC(h_3, 18)=23  $.

When \(dt\) increases by $ \Delta t $,  $dt=1  $. At this moment,  $ T_{h_2}(1)=2 $ and $ tc(h_1h_2h_3)=25 $. When $ dt=2 $, $ T_{h_2}(2)=3 $, and $ tc(h_1h_2h_3)=21 $. Furthermore, when  $ dt=10 $, $ T_{h_2}(10)=25 $ and $tc(h_1h_2h_3)=115  $.

From the above computations, it is evident that when  \(dt\) is 0, 1, 2, and 10 respectively, \(T_{h2}(dt)\) is 3, 2, 3, and 25, respectively. The route cost, i.e., \(tc(h_1h_2h_3)\), are 23, 25, 21, 115 respectively. It can be observed  that as \(dt\) increases from 0 to 1, \(T_{h_2}(dt)\) decreases from 3 to 2. At the same time, as \(dt\) increases continually from 0, the route cost initially  increases, then decreases, and subsequently  increases again, which demonstrates a non-unimodal relationship, as shown in Fig. \ref{fig_m}(c). 
\begin{table*}[!htb]	
	\scriptsize
	\caption{Results on the $gdb$ benchmark test set in the 2LP in terms of costs of solutions. For each instance, the average performance of an algorithm is indicated in bold if it is the best among all comparison algorithms. ``$\ddagger$" and ``$\dagger$" indicate MAENS-GN performs significantly better than and equivalently to the corresponding algorithm based on 20 independent runs according to the Wilcoxon rank-sum test with significant level $ p=0.05 $, respectively. }
	\label{tab:1}	
	\centering
	\noindent \makebox[\textwidth][c]
	{
		\begin{tabular}{ccccccccccc}
			\hline  \noalign{\smallskip} 
			Instances &$|R|$ & LB  & \multicolumn{3}{c}{MAENS-GN}  & \multicolumn{3}{c}{VND}  &\multicolumn{2}{c}{VND*}          \\
			 \cmidrule(lr){4-6} \cmidrule(lr){7-9} \cmidrule(lr){10-11}
			\noalign{\smallskip} 
			&        &     & Ave(std)        & Best  & Time    & Ave(std)         & Best & Time              & Ave(std)         & Best     \\
			\hline  \noalign{\smallskip}
			gdb1   &22 & 316 & \textbf{316.0(0.00)}     & 316   & 0.799   & 363.0(0.00)$\ddagger$     & 363  & 0.102    &359.5(5.82)$\ddagger$ 	 &340\\
			\noalign{\smallskip}
			gdb2   &26 & 339 & \textbf{345.2(0.65)}     & 345   & 8.505   & 379.6(2.82)$\ddagger$      & 375  & 0.149     &372.4(9.30)$\ddagger$ 	&346\\
			\noalign{\smallskip}
			gdb3  &22  & 275 & \textbf{289.0(0.00) }    & 289   & 6.131   & 315.4(0.49)$\ddagger$      & 315  & 0.088     &307.7(8.25)$\ddagger$ 	&289\\
			\noalign{\smallskip}
			gdb4   &19 & 287 & \textbf{287.0(0.00) }    & 287   & 1.118   & 342.0(0.00)$\ddagger$       & 342  & 0.056     &334.9(7.65)$\ddagger$ 	&317\\
			\noalign{\smallskip}
			gdb5  &26  & 377 & \textbf{381.7(3.84)}     & 377   & 7.287   & 417.0(0.00)$\ddagger$      & 417  & 0.121     &407.4(9.83)$\ddagger$ 	&388\\
			\noalign{\smallskip}
			gdb6  &22  & 298 & \textbf{299.3(5.67)}     & 298   & 1.395   & 344.0(0.00)$\ddagger$       & 344  & 0.085     &339.2(7.65)$\ddagger$ 	&324\\
			\noalign{\smallskip}
			gdb7  &22  & 325 & \textbf{325.0(0.00)}     & 325   & 1.548   & 377.8(4.28)$\ddagger$       & 364  & 0.088     &359.0(13.57)$\ddagger$ 	&335\\
			\noalign{\smallskip}
			gdb8  &46  & 348 & \textbf{357.1(1.48)}     & 355   & 34.872  & 395.0(0.00)$\ddagger$       & 395  & 0.588    &389.8(5.78)$\ddagger$ 	&370\\
			\noalign{\smallskip}
			gdb9  &51  & 303 & \textbf{312.0(3.19)}     & 309   & 38.796  & 358.2(4.59)$\ddagger$      & 346  & 0.836     &343.2(4.36)$\ddagger$ 	&336\\
			\noalign{\smallskip}
			gdb10 &25  & 275 & \textbf{278.6(4.04) }    & 275   & 5.327   & 313.8(5.08)$\ddagger$      & 309  & 0.117     &306.6(7.45)$\ddagger$ 	&290\\
			\noalign{\smallskip}
			gdb11 &45  & 395 & \textbf{403.8(2.14)}     & 399   & 32.816  & 446.1(2.86)$\ddagger$       & 435  & 0.671     &429.8(4.36)$\ddagger$ 	&421\\
			\noalign{\smallskip}
			gdb12 &23  & 458 & \textbf{458.0(0.00)}     & 458   & 1.061   & 503.0(0.00)$\ddagger$      & 503  & 0.082     &503.0(0.00)$\ddagger$     &503\\
			\noalign{\smallskip}
			gdb13 &28  & 536 & \textbf{545.4(4.75) }    & 536   & 9.129   & 586.0(6.49)$\ddagger$      & 578  & 0.150     &569.7(4.52)$\ddagger$ 	&560\\
			\noalign{\smallskip}
			gdb14 &21  & 100 & \textbf{101.5(1.50)}     & 100   & 4.179   & 109.0(0.00)$\ddagger$      & 109  & 0.072     &107.0(1.73)$\ddagger$ 	&104\\
			\noalign{\smallskip}
			gdb15 &21  & 58  & \textbf{58.0(0.00)}      & 58    & 0.255   & 62.0(0.00)$\ddagger$       & 62   & 0.072     &61.5(1.07)$\ddagger$ 	&58\\
			\noalign{\smallskip}
			gdb16 &28  & 127 & \textbf{129.0(0.00)}     & 129   & 10.749  & 137.0(0.00)$\ddagger$      & 137  & 0.146    &133.1(1.04)$\ddagger$ 	&131\\
			\noalign{\smallskip}
			gdb17 &28  & 91  & \textbf{91.0(0.00)}      & 91    & 0.123   & 92.4(0.92)$\ddagger$       & 91   & 0.161    &92.8(0.60)$\ddagger$ 	& 91\\
			\noalign{\smallskip}
			gdb18 &36  & 164 & \textbf{169.4(1.56)}     & 164   & 14.403  & 178.2(3.18)$\ddagger$      & 172  & 0.307     &173.2(1.44)$\ddagger$  &171\\
			\noalign{\smallskip}
			gdb19 &11  & 55  & \textbf{55.0(0.00)}      & 55    & 0.053   & 63.0(0.00)$\ddagger$      & 63   & 0.020      &62.1(1.84)$\ddagger$ 	&55\\
			\noalign{\smallskip}
			gdb20 &22  & 121 & \textbf{122.8(0.60)}     & 121   & 6.309   & 125.0(0.00)$\ddagger$      & 125  & 0.070      &124.8(0.51)$\ddagger$ 	&123\\
			\noalign{\smallskip}
			gdb21 &33  & 156 & \textbf{157.9(0.54)}     & 156   & 12.504  & 169.2(2.17)$\ddagger$      & 166  & 0.236    &165.4(2.15)$\ddagger$ 	&162\\
			\noalign{\smallskip}
			gdb22 &44  & 200 & \textbf{202.2(0.60)}     & 202   & 25.411  & 206.0(0.00)$\ddagger$      & 206  & 0.538     &205.2(1.01)$\ddagger$ 	&203\\
			\noalign{\smallskip}
			gdb23 &55  & 233 & \textbf{237.5(0.92)}     & 236   & 41.476  & 244.0(0.00)$\ddagger$      & 244  & 1.003     &243.8(0.36)$\ddagger$ 	&243\\
			\noalign{\smallskip}
			w-d-l  &-   & -   & -       & -     & -       & 23-0-0                      & -    & -         & 23-0-0	   &-\\		
			\noalign{\smallskip}
			No.best &-  & -   & 9             & 23    & -       & 0                           & 1    & -         &0	       &4\\
			\noalign{\smallskip}
			Ave.PDR &-  & -   & 1.37\%        & -     & -       & 10.82\%                     & -    & -         &8.74\% 	   &- \\
			\noalign{\smallskip}
			Ave.Time &-  & -   & -             & -     & 11.489  & -                           & -    & 0.25      &-          & -\\
			\hline  \noalign{\smallskip}
		\end{tabular}
	}
\end{table*}

\subsubsection{Methods for optimizing  the vehicle departure time}
\label{solvtime}
In the 2LP, time 0 is the optimal vehicle departure time. 

In the 3LP, when the slope $ |k|\leq 1 $, the route cost and the departure time of the vehicle exhibit  a (similar) unimodal or non-decreasing relationship. In this case,  GSS \cite{OVERHOLT1967FIBONACCI} is used to  quickly identify  the optimal vehicle departure time for each route in the routing plan. GSS,  characterized as an interval search method,  initiates with the entire planning horizon, i.e., [0, \(T\)], and successively narrows down the search interval using the golden ratio until the length of the search interval becomes smaller than a threshold $ \varepsilon $. At this point, the optimal vehicle departure time is determined. In addition, GSS is an algorithm with low time complexity. For more details on GSS, please refer to \cite{OVERHOLT1967FIBONACCI}. \\
\indent In the 3LP, when the slope $|k|>1$ , the route cost and the departure time of the vehicle exhibit a non-unimodal relationship, and GSS is no longer applicable. In such a scenario, an evolutionary algorithm, that is, NCS \cite{tang2016negatively}, is employed  to obtain  the (approximately) optimal vehicle departure time. NCS is good at dealing with complex optimization problems. It mainly maintains the search process of multiple individuals in parallel, and models the individual search process as a probability distribution. For more details on NCS, please refer to \cite{tang2016negatively}.
\begin{table*}[!htb]	
	\caption{Results on the $egl$ benchmark test set in the 2LP in terms of costs of solutions.  For each instance, the average performance of an algorithm is indicated in bold if it is the best among all comparison algorithms. ``$\ddagger$" and ``$\dagger$" indicate MAENS-GN performs significantly better than and equivalently to the corresponding algorithm based on 20 independent runs according to the Wilcoxon rank-sum test with significant level $p=0.05$, respectively.}
	\label{tab:2}
	\scriptsize
	\centering
	
	\noindent \makebox[\textwidth][c]
	{
		\begin{tabular}{ccccccccccc}
			\hline  \noalign{\smallskip}
			Instances &$|R|$ & LB    & \multicolumn{3}{c}{MAENS-GN} & \multicolumn{3}{c}{VND}           & \multicolumn{2}{c}{VND*}  \\
			\cmidrule(lr){4-6} \cmidrule(lr){7-9}  \cmidrule(lr){10-11}
			&     &           & Ave(std)        & Best   & Time         & Ave(std)                   & Best   & Time    & Ave(std)              & Best         \\
			\hline \noalign{\smallskip}
			egl-e1-A &51& 3548  & \textbf{3550.6(8.97)}     & 3548  & 28.973       & 4323.8(110.19)$\ddagger$    & 4018  & 0.687   & 4033.2(94.85)$\ddagger$		& 3900\\
			\noalign{\smallskip}
			egl-e1-B &51& 4498  & \textbf{4576.9(18.31)}    & 4557  & 74.811       & 5006.6(35.66)$\ddagger$     & 4877  & 0.902   & 4760.1(80.10)$\ddagger$		& 4573\\
			\noalign{\smallskip}
			egl-e1-C &51& 5595  & \textbf{5677.6(59.91)}    & 5595  & 66.548       & 6178.4(71.04)$\ddagger$     & 5923  & 1.017   & 6091.4(84.25)$\ddagger$		& 5907\\
			\noalign{\smallskip}
			egl-e2-A &72& 5018  & \textbf{5102.4(64.53)}    & 5018  & 189.315      & 5886.0(158.29)$\ddagger$    & 5558  & 2.552   & 5425.3(73.00)$\ddagger$		& 5307\\
			\noalign{\smallskip}
			egl-e2-B &72& 6317  & \textbf{6424.9(26.91)}   & 6355  & 163.441       & 7114.4(75.27)$\ddagger$     & 6869  & 1.765   & 6703.0(66.97)$\ddagger$		& 6554\\
			\noalign{\smallskip}
			egl-e2-C &72& 8335  & \textbf{8599.0(64.46)}    & 8474  & 163.261      & 8910.6(5.88)$\ddagger$      & 8885  & 1.768   & 8779.8(101.86)$\ddagger$		& 8548\\
			\noalign{\smallskip}
			egl-e3-A &87& 5898  & \textbf{6099.4(77.63)}    & 6002  & 335.078      & 6483.0(0.00)$\ddagger$      & 6483  & 3.104   & 6475.2(34.00)$\ddagger$		& 6327\\
			\noalign{\smallskip}
			egl-e3-B &87& 7744  & \textbf{8040.6(71.43)}    & 7909  & 251.941      & 8819.0(0.00)$\ddagger$      & 8819  & 4.752   & 8491.9(94.51)$\ddagger$		& 8264\\
			\noalign{\smallskip}
			egl-e3-C &87& 10244 & \textbf{10460.6(54.10)}   & 10390 & 213.761      & 11485.4(185.13)$\ddagger$   & 10865 & 2.642   & 10928.1(141.76)$\ddagger$		& 10429\\
			\noalign{\smallskip}
			egl-e4-A &98& 6408  & \textbf{6710.6(100.79)}   & 6580  & 364.345      & 7526.5(19.66)$\ddagger$     & 7460  & 6.271   & 7222.0(81.88)$\ddagger$		& 7075\\
			\noalign{\smallskip}
			egl-e4-B &98& 8935  & \textbf{9329.9(78.16)}   & 9165  & 318.551       & 10251.3(170.85)$\ddagger$   & 9876  & 5.741   & 9748.6(101.41)$\ddagger$		& 9438\\
			\noalign{\smallskip}
			egl-e4-C &98& 11512 & \textbf{11992.6(91.87)}   & 11854 & 264.948      & 12614.8(14.17)$\ddagger$    & 12553 & 5.560    & 12484.6(114.52)$\ddagger$		& 12225\\
			\noalign{\smallskip}
			egl-s1-A &75& 5018  & \textbf{5145.9(102.02)}   & 5018  & 218.771      & 5804.3(7.41)$\ddagger$      & 5772  & 2.148   & 5613.8(154.15)$\ddagger$		& 5267\\
			\noalign{\smallskip}
			egl-s1-B &75& 6388  & \textbf{6576.6(106.89)}  & 6394  & 187.171       & 6856.0(0.00)$\ddagger$      & 6856  & 2.677   & 6856.0(0.00)$\ddagger$	        & 6856\\
			\noalign{\smallskip}
			egl-s1-C &75& 8518  & \textbf{8783.8(62.16)}    & 8677  & 162.748      & 9643.8(17.13)$\ddagger$     & 9573  & 2.794   & 9479.2(96.17)$\ddagger$		& 9256\\
			\noalign{\smallskip}
			egl-s2-A &147& 9825  & \textbf{10490.1(60.67)}   & 10392 & 1172.03      & 11414.2(86.52)$\ddagger$    & 11037 & 23.580   & 11204.4(109.18)$\ddagger$		& 11015\\
			\noalign{\smallskip}
			egl-s2-B &147& 13017 & \textbf{13887.2(119.45)} & 13698 & 795.669       & 14680.0(0.00)$\ddagger$     & 14680 & 14.923  & 14661.8(30.65)$\ddagger$		& 14567\\
			\noalign{\smallskip}
			egl-s2-C &147& 16425 & \textbf{17392.6(139.47)}  & 17165 & 761.257      & 18408.6(14.82)$\ddagger$    & 18344 & 18.951  & 18327.3(83.45)$\ddagger$		& 18174\\
			\noalign{\smallskip}
			egl-s3-A &159& 10165 & \textbf{10733.0(99.71)}     & 10593 & 1417.264   & 11764.7(75.41)$\ddagger$    & 11436 & 19.466  & 11406.1(100.94)$\ddagger$		& 11207\\
			\noalign{\smallskip}
			egl-s3-B &159& 13648 & \textbf{14418.0(125.53)}    & 14168 & 1046.018   & 15858.4(169.85)$\ddagger$   & 15367 & 27.001  & 15325.0(93.46)$\ddagger$		& 15111\\
			\noalign{\smallskip}
			egl-s3-C &159& 17188 & \textbf{18269.3(151.44)}  & 17983 & 902.040       & 19666.0(0.00)$\ddagger$     & 19666 & 22.921  & 19277.8(155.11)$\ddagger$		& 19029\\
			\noalign{\smallskip}
			egl-s4-A &190& 12153 & \textbf{13197.1(152.21)}  & 12932 & 1531.802     & 14918.2(319.64)$\ddagger$   & 14285 & 46.072  & 14142.2(139.63)$\ddagger$	    & 13762\\
			\noalign{\smallskip}
			egl-s4-B &190& 16113 & \textbf{17260.1(91.44)}   & 17121 & 1499.895     & 18407.2(12.20)$\ddagger$    & 18354 & 34.835  & 18242.6(144.89)$\ddagger$		& 17943\\
			\noalign{\smallskip}
			egl-s4-C &190& 20430 & \textbf{21796.0(188.09)}    & 21427 & 1348.113   & 22996.4(2.83)$\ddagger$     & 22984 & 37.787  & 22945.9(113.44)$\ddagger$		& 22513\\
			\noalign{\smallskip}
			w-d-l    &-  & -     & -           & -      & -         & 24-0-0                      & -     & -       & 24-0-0  &-\\			
			\noalign{\smallskip}
			No.best  &-  & -     & 0                & 24     & -         & 0                           & 0     & -       &0        &0\\
			\noalign{\smallskip}
			Ave.PDR  &-  & -     & 4.18\%          & -      & -          & 13.80\%                     & -     & -       &10.48\% &-\\
			\noalign{\smallskip}
			Ave.Time &-  & -     & -                & -     & 561.573    & -                           & -     & 12.08   &-       &-\\
			\hline
		\end{tabular}
	}
\end{table*}

\subsection{Novel initialization strategy}
\label{indinit}

\begin{table*}[!htb]	
	\caption{Results on the $Solomon$-25 benchmark test set in the 3LP in terms of costs of solutions. For each instance, the average performance of an algorithm is indicated in bold if it is the best among all comparison algorithms (except CG). ``$\ddagger$" and ``$\dagger$" indicate MAENS-GN performs significantly better than and equivalently to the corresponding algorithm based on 20 independent runs according to the Wilcoxon rank-sum test with significant level $p=0.05  $, respectively.}
	\label{tab:3}
	\scriptsize
	\centering
	
	\noindent \makebox[\textwidth][c]
	{
		\begin{tabular}{cccccccccc}
			\hline  \noalign{\smallskip}
			Instances            & \multicolumn{2}{c}{CG}                & \multicolumn{3}{c}{MAENS-GN}                                                       &                       \multicolumn{3}{c}{VND-GN}                                                      \\
			\cmidrule(lr){2-3}  \cmidrule(lr){4-6}  \cmidrule(lr){7-9} 
			\noalign{\smallskip}
			                    & Cost                 & Time                 & Ave(std)                  & Best                 & Time                 & Ave(std)                          & Best                 & Time                     \\
			\hline   \noalign{\smallskip}
			r101                 & 893.7                & 258.4                & \textbf{1014.3(10.86)}              &  1003.4                & 10.725               & 1089.9(34.49)$\ddagger$              & 1050.3                & 0.100                   \\                                         
			\noalign{\smallskip}
			r102                 & 806.0                & 410.0                & \textbf{865.5(11.79)}             & 853.8                & 9.955                & 908.1(3.68)$\ddagger$              & 900.8               & 0.095                    \\
                                                                                        						
			\noalign{\smallskip}
			r103                 & 726.1                & 80.5                 &\textbf{ 770.7(2.14) }             & 767.0                & 10.525                & 819.2(6.22)$\ddagger$               & 806.8                & 0.095                      \\                                                                                   		               	                    	                 		                         	
			\noalign{\smallskip}
			r104                 & 690.0                & 571.1                & \textbf{715.2(8.01)}              & 703.0                & 6.865                 &  767.8(27.37)$\ddagger$              & 736.8                & 0.086                       \\                                            
			\noalign{\smallskip}
			r105                 & 780.4                & 118.2                & \textbf{839.7(3.18)}              &  836.3                & 10.140               & 946.1(38.53)$\ddagger$              & 853.9                & 0.090                        \\                                       
			\noalign{\smallskip}
			r106                 & 721.4                & 236.9                & \textbf{744.0(4.61) }             & 740.1	                & 9.595               & 864.6(7.33)$\ddagger$              & 847.3                & 0.097                       \\                                        

			\noalign{\smallskip}
			r107                 & 679.5                & 137.8                & \textbf{700.5(2.81)}              &  697.5               & 9.850               & 760.0(9.95)$\ddagger$               & 740.1                  & 0.092                         \\                                          
			\noalign{\smallskip}
			r108                 & 647.7                & 232.9                & \textbf{663.1(10.17)  }           & 653.1                  & 6.900                & 712.0(17.36)$\ddagger$              & 689.4               & 0.087                         \\                                       	
			\noalign{\smallskip}
			r109                 & 691.3                & 71.7                 & \textbf{744.0(2.33) }             & 734.7                & 9.460               & 822.5(22.36)$\ddagger$               & 777.8                & 0.087                     \\                                                                                        		          	                	                 		                   	

			\noalign{\smallskip}
			r110                 & 668.8                & 125.9                & \textbf{716.8(10.05)}              & 697.2                & 15.925               & 758.8(23.82)$\ddagger$                       &  729.8                & 0.093                       \\                                                                               		      	                    	                      		                 	

			\noalign{\smallskip}
			r111                 & 676.3                & 157.8                & \textbf{710.5(5.16)}              & 707.9                &  16.575                 & 778.9(16.80)$\ddagger$                 & 746.3                & 0.094                     \\                                                                                    		             	                    	                  			

			\noalign{\smallskip}
			r112                 & 643.0                & 356.1                & \textbf{666.3(5.17)}              & 656.0                & 15.925               & 718.3(13.18)$\ddagger$               & 696.1                  & 0.088                    \\                                                                                 		                  	                	                        			

			\noalign{\smallskip}
			rc101                & 623.8                & 62.4                 & \textbf{648.4(3.07) }             & 643.9                & 19.290                & 761.3(58.66)$\ddagger$                 & 681.7                & 0.090                      \\                                                                                          		                 		                               			

			\noalign{\smallskip}
			rc102                & 598.3                & 179.8                & \textbf{611.5(5.07)}              & 606.6                & 20.185               & 803.9(36.15)$\ddagger$                & 761.5                & 0.103                      \\                                                                                        				                                          			

			\noalign{\smallskip}
			rc103                & 585.1                & 941.5                & \textbf{607.8(4.96)}              & 598.1                & 20.005                & 738.2(26.23)$\ddagger$                 & 645.4                & 0.089                      \\                                                                                         				                                            			

			\noalign{\smallskip}
			rc104                & 572.4                & 3183.1               & \textbf{603.3(3.52)}              & 595.4                & 20.705               & 609.0(0.00)$\ddagger$                & 609.0                & 0.077                     \\                                                                                        				                                          			

			\noalign{\smallskip}
			rc105                & 623.1                & 743.7                & \textbf{656.3(9.24)}              & 641.1                & 21.335               & 813.6(62.44)$\ddagger$                 & 677.3                & 0.082                       \\                                                                                    		                		                              			
			\noalign{\smallskip}
			rc106                & 588.7                & 372.6                & \textbf{605.6(0.00)}              & 605.6                & 18.370                & 631.3(0.00)$\ddagger$                 & 631.3                & 0.084                         \\                                                                                         		        		                                 			

			\noalign{\smallskip}
			rc107                & 548.3                & 824.8                & \textbf{567.0(2.15)}              & 565.4               & 12.530                & 628.5(11.66)$\ddagger$                  & 600.1                & 0.079                        \\                                                                                       				                                             			

			\noalign{\smallskip}
			rc108                & 544.5                & 3457.1               & \textbf{553.5(4.42)}              & 546.8                & 12.630                & 579.7(13.99)$\ddagger$                   & 555.9                & 0.075                        \\                                                                                  				                                              			

			\noalign{\smallskip}
			w-d-l                & 0-0-20                   & -                    & -                  & -                    & -                    & 20-0-0                        & -                    & -                      \\
			\noalign{\smallskip}
			No.best              & 20                   & -                    & 0                      & 0                    & -                    & 0                              & 0                    & -                      \\
			\noalign{\smallskip}
			Ave.PDR              & 0.00\%                  & -                    & 4.95\%                 & -                    & -                    & 16.44\%                         & -                    & -                        \\
			\noalign{\smallskip}
			Ave.Time             & -                    & 626.115              & -                     & -                    & 13.875              & -                               & -                    & 0.089                   \\
			\hline  \noalign{\smallskip}
		\end{tabular}
	}
\end{table*}

When constructing a route in the individual initialization strategy of MAENS, the nearest task to the end of the route  needs to be selected to add to the route. If more than one task satisfies the condition simultaneously, these tasks form a task set termed  $Tasks_{nearst}$.
At this time, the tie-breaking rule selected is to randomly select one from $Tasks_{nearst}$. 
However, this stochastic tie-breaking rule is a little blind in the CARPTDSC. Because in the CARPTDSC, the service cost of the task is time-dependent, meaning that  the task's service cost varies based on its position (i.e., the time of beginning of service) in the route. Consequently, this stochastic tie-breaking rule may select a task with the larger service cost caused by its time of beginning of service  being distant  from its optimal service time interval, such as the time of beginning of service of \(h_1\) is 10 in Fig. \ref{fig4}(a).\\
\indent In this context, we  design a specialized tie-breaking rule that takes into account the time-dependent property of CARPTDSC  in the individual initialization strategy.  If choosing the task with the smallest service cost from  $Tasks_{nearst}$, the time of beginning of service of the selected task may fall within or not far from  its optimal service time interval (e.g., the time interval \([bt, et]\) in Fig. \ref{fig0}(a)). In this way, it is avoided that the task is in another position (i.e., the time of beginning of service  of the task is far away from its  optimal service time interval) of the route, resulting in a larger service cost. Because each task needs to be inserted into the routing plan and only once, selecting the task being in a good position can avoid the increase of service cost caused by its being in a bad position. In this way,  the further increase of the route cost is avoided. However, if doing this, the individuals produced in the  initial population will all be the same, meaning no diversity. So this is not a good method. To address this,  we devised a roulette wheel selection-based \cite{holland1975adaptation} approach where the task with the smaller service cost will have a greater probability of being selected. Therefore, this method not only considers diversity but also makes tasks with better quality more likely to be selected.\\
\indent  Here we will select a task from  $Tasks_{nearst} $ by the roulette wheel selection strategy and   one score is associated with each task in $Tasks_{nearst} $. The higher the score of the task, the higher the probability of being selected. Due to the need to make the
service cost of the task lower, the probability of being selected is higher. Consequently, we define  the score of task $x $, denoted as $score(x)  $, to be equal to the reciprocal of the task's time-dependent service cost, as expressed  in Eq. (\ref{equ:9}).
\begin{equation} \label{equ:9}
	score(x)=\frac{1}{SC(x, t)},
\end{equation}
where $ SC(x, t) $ is the time-dependent service cost of task $x$. The probability of the task $x_i$ being selected, denoted as $p(x_i)$, is the proportion of the score of task $x_i$ to the sum of the scores of all tasks in $tasks_{nearst}$, as shown in Eq. (\ref{equ:10}).
\begin{equation}\label{equ:10}
p(x_i)=\frac{score(x_i)}{\sum_{j=1}^{N} score(x_j)},
\end{equation}
where $ N $ is the size of $ Tasks_{nearst} $. The detailed process of the roulette wheel selection will not be introduced in detail below. For more details on roulette wheel selection,  please refer to the \cite{holland1975adaptation}.
Below is an illustrative  example.\\
\indent Assume that there are only two tasks \(h_2\) and \(h_3\) (in Fig. \ref{fig4}) that have not been inserted into the initial routing plan \(X_{init}\), and both tasks are in \(Tasks_{nearst}\), and the two tasks are connected directly. Therefore, the order of inserting into \(X_{init}\) of these two tasks  depends on the task selected firstly from \(Tasks_{nearst}\).  Because the arrival time of the tasks in the \(Tasks_{nearst}\) is the same, assume that it is 12. At this time, the service costs of \(h_2\) and \(h_3\) can be obtained  to be \(SC(h_2, 12)=1\) and \(SC(h_3,12)=5\), respectively. If roulette wheel selection is used, the probabilities of \(h_2\) and \(h_3\) being selected is \((1/1)/(1/1+1/5)\approx0.83\) and 0.17, respectively.  At this time, if \(h_2\) is inserted firstly and then \(h_3\), then \(SC(h_2, 12)=1\) and \(SC(h_3, 13)=3\). The total cost of the  sequence of tasks \(h_2h_3\)  is \(tc(h_2h_3)=1+3=4\).  If \(h_3\) is inserted firstly and then \(h_2\), then \(SC(h_3, 12)=5\) and \(SC(h_2, 17)=11\), leading to that a total cost of \( h_3h_2 \) is \(tc(h_3h_2)=16\). Using the roulette wheel selection, the probability of the total cost of the sequence \(h_2h_3\) being 4 is 0.83, and the probability of the total cost of \(h_3h_2\) being 16 is 0.17. In contrast, if  the stochastic selection  method is used, the probabilities of  the total costs of the  sequence of tasks being 4 and 16 would both be 0.5. \\
\indent Therefore, the roulette wheel selection exhibits a higher likelihood of acquiring  an initial solution with high quality.
\section{Experimental studies}
\label{expe}
The experiments aim to verify whether MAENS-GN has effectively  addressed the two research questions (RQ1 and RQ2). Consequently, the verification of the performance of MAENS-GN is divided into two aspects.
\begin{enumerate}
	\item In the 2LP, the quality of approximate solutions obtained by existing algorithms is not high, but does MAENS-GN significantly improve the quality of approximate solutions? 
	\item In the 3LP, the existing exact algorithms can only handle small-scale instances (number of tasks \(\leq\)40), whereas what is the quality of the solutions of MAENS-GN on larger-scale instances in addition to small-scale instances?
\end{enumerate}
Additionally, we also assessed  the contribution of the components in MAENS-GN, including  GSS, NCS and  the novel initialization strategy. In each instance, each algorithm was run independently 20 times to ensure robust evaluations.
\begin{table*}[!htb]	
	\caption{Results on the $Solomon$-35 benchmark test set in the 3LP in terms of costs of solutions.  For each instance, the average performance of an algorithm is indicated in bold if it is the best among all comparison algorithms (except CG). ``$\ddagger$" and ``$\dagger$" indicate MAENS-GN performs significantly better than and equivalently to the corresponding algorithm based on 20 independent runs according to the Wilcoxon rank-sum test with significant level $p=0.05  $, respectively.}
	\label{tab:4}
	\scriptsize
	\centering
	
	\noindent \makebox[\textwidth][c]
	{
		\begin{tabular}{ccccccccc}
			\hline  \noalign{\smallskip}
			Instances            & \multicolumn{2}{c}{CG}                      & \multicolumn{3}{c}{MAENS-GN}                                                       & \multicolumn{3}{c}{VND-GN}                                                         \\
				\cmidrule(lr){2-3}  \cmidrule(lr){4-6}  \cmidrule(lr){7-9} 
			\noalign{\smallskip}
		   	                     & Cost                 & Time                 & Ave(std)                 & Best                 & Time                 & Ave(std)                  & Best                 & Time                                     \\
			\hline  \noalign{\smallskip}
			r101                 & 1168.3               & 4165.8               & \textbf{1316.9(12.38)}           & 1293.9               & 42.405               & 1373.4(6.25)$\ddagger$                & 1369.5               & 0.312                            \\
                                                                                 				                                           			
			\noalign{\smallskip}
			r102                 & 1048.0               & 6574.1               & \textbf{1097.3(7.34)}          & 1085.0               & 45.620                & 1265.9(44.77)$\ddagger$               & 1158.8               & 0.282                            \\                                                                                  		            	            	                  			
			\noalign{\smallskip}
			r103                 & 929.2                & 5129.9               & \textbf{981.1(8.18)}            & 969.4                & 43.460                & 1126.1(34.16)$\ddagger$              & 1071.4               & 0.227                  \\
                                                                              		      	                   	                			
			\noalign{\smallskip}
			r104                 & 837.8                & 6210.9               & \textbf{862.4(7.97)}           &  853.4                & 49.230                & 921.9(9.63)$\ddagger$                & 898.3                  & 0.250                            \\
                                                                               		               		                            			
			\noalign{\smallskip}
			r105                 & 1017.1               & 2926.2               & \textbf{ 1108.2(10.03)}           & 1091.9               & 46.450                & 1186.5(17.82)$\ddagger$              & 1156.7               & 0.300                   \\
                                                                               		        		                                			
			\noalign{\smallskip}
			r106                 & 930.8                & 2247.1               & \textbf{968.5(9.10)}            & 956.4                & 43.235               & 1090.6(27.85)$\ddagger$                & 1045.2               & 0.277                          \\
                                                                                				                                               		
			\noalign{\smallskip}
			r107                 & 874.8                & 8957.4               & \textbf{911.8(10.55)}           & 896.2                & 44.770                & 1011.2(21.72)$\ddagger$                &944.1                & 0.260                           \\
                                                                                				                                      			
			\noalign{\smallskip}
			r108                 & 805.7                & 14698.5              & \textbf{826.5(5.89)}            & 820.3                & 47.030                & 887.7(28.64)$\ddagger$                 & 841.0                & 0.259                          \\
                                                                                 				                                        			
			\noalign{\smallskip}
			r109                 & 914.8                & 8044.2               & \textbf{951.7(9.51) }          & 942.0                & 43.595                 & 1010.5(14.26)$\ddagger$                 & 976.1                & 0.274                           \\
                                                                                 			                                        			
			\noalign{\smallskip}
			r110                 & 893.6                & 16259.7              & \textbf{910.5(8.71)}           & 902.2                & 33.580               & 993.1(40.58)$\ddagger$                    & 948.5                & 0.257                           \\
                                                                                		         		                                		                      	
			\noalign{\smallskip}
			r111                 & 864.4                & 9957.6               & \textbf{895.6(5.86)}            & 887.6                & 35.015                & 913.8(1.96)$\ddagger$                 & 907.2                  & 0.290                             \\
                                                                               		          		                               			
			\noalign{\smallskip}
			r112                 & 834.3                & 10377.9              & \textbf{839.6(2.84)}           & 835.1                & 35.230                & 887.9(8.24)$\ddagger$                     & 854.3                & 0.239                           \\
                                                                                    		     		                                    			
			\noalign{\smallskip}
			w-d-l                & 0-0-12                    & -                    & -                  & -                    & -                    & 12-0-0                           & -                    & -                                    \\
			
			\noalign{\smallskip}
			No.best              & 12                   & -                    & 0                        & 0                    & -                    & 0                                 & 0                    & -                       \\
			\noalign{\smallskip}
			Ave.PDR              & 0.00\%                  & -                    & 4.66\%                   & -                    & -                    & 13.58\%                          & -                    & -                        \\
			
			\noalign{\smallskip}
			Ave.Time             & -                    & 7962.442             & -                      & -                    & 42.468             & -                                 & -                 & 0.269                     \\
			\hline  \noalign{\smallskip}
		\end{tabular}
	}
\end{table*}

\subsection{Experimental setup}
\label{setup}
Two distinct  types of benchmark test sets, i.e., test sets in the 2LP and  3LP,  were utilized for evaluations. The $gdb$  and $egl$ in the 2LP, which already exist in \cite{tagmouti2010variable}, were modified based on the static test sets $gdb  $ \cite{dearmon1981comparison} and $ egl $ \cite{eglese1994routeing}, i.e., two-segment linear  functions were added. In the 3LP, Solomon's test set \cite{solomon1987algorithms}, as used in \cite{tagmouti2007arc}, was adopted. In \cite{tagmouti2007arc}, 25-customer, 35-customer, and 40-customer instances in Solomon's test set were used, which are denoted as $ Solomon $-25, $ Solomon $-35, $ Solomon $-40 here. In addition to Solomon's test set, the same method in \cite{tagmouti2010variable} was used to apply the three-segment linear function on the static benchmark test set $ gdb $ and $ egl $, making them  time-dependent test sets in the 3LP. In addition, a real-world  $jd  $ dataset was also used, which is derived from the Jingdong logistics distribution system. The $ jd $ set is similar to  Solomom's test set,  being also VRP test set with time windows. For convenience, the absolute values of the slopes in a instance of 3LP  are the same, such as 2 and -2 in Fig. \ref{fig4}(a), as  done in the \cite{tagmouti2007arc}. 
The slope in the test sets in the 2LP is  1 \cite{tagmouti2010variable}. The slopes in Solomon's test set are 1 and -1 \cite{tagmouti2007arc}.  Furthermore, the absolute values of the slopes in the instances of other test sets, i.e., \(gdb\), \( egl\), and \(jd\) in the 3LP, were randomly chosen from the set \{0.3, 0.5, 1, 2, 3\}. 
\begin{table*}[!htb]	
	\caption{Results on the $Solomon$-40 benchmark test set in the 3LP in terms of costs of solutions. For each instance, the average performance of an algorithm is indicated in bold if it is the best among all comparison algorithms (except CG). ``$\ddagger$" and ``$\dagger$" indicate MAENS-GN performs significantly better than and equivalently to the corresponding algorithm based on 20 independent runs according to the Wilcoxon rank-sum test with significant level $p=0.05  $, respectively.}
	\label{tab:5}
	\scriptsize
	\centering	
	\noindent \makebox[\textwidth][c]
	{
		\begin{tabular}{ccccccccc}
			\hline  \noalign{\smallskip}
			Instances            & \multicolumn{2}{c}{CG}                      & \multicolumn{3}{c}{MAENS-GN}                                                       & \multicolumn{3}{c}{VND-GN}                                                                     \\
				\cmidrule(lr){2-3}  \cmidrule(lr){4-6}  \cmidrule(lr){7-9} 
			                     & Cost                 & Time                 & Ave(std)                 & Best                 & Time                 & Ave(std)                          & Best                 & Time                                      \\
			\hline  \noalign{\smallskip}
			r101                 & 1318.4               & 4467.3               & \textbf{1479.8(14.61) }          & 1449.0               & 55.340                 & 1602.5(23.20)$\ddagger$               & 1531.0               & 0.462                            \\
             
			\noalign{\smallskip}
			r102                 & 1194.8               & 18641.9              & \textbf{1268.1(11.38) }          & 1255.1               & 59.620               & 1403.6(37.67)$\ddagger$               & 1338.2               & 0.399                           \\                                                                               							
 
			\noalign{\smallskip}
			r103                 & -                    & -                    & \textbf{1113.8(13.83)}           & 1099.9               & 59.500                & 1238.6(18.60)$\ddagger$               & 1189.6               & 0.385                            \\                                                                               							
           
			\noalign{\smallskip}
			r104                 & -                    & -                    &\textbf{975.0(4.40)}           & 970.2                & 63.320               & 1065.7(15.82)$\ddagger$                & 1027.3               & 0.325                              \\                                                                               							
    
			\noalign{\smallskip}
			r105                 & 1171.2               & 3780.2               & \textbf{1267.7(7.55)}          & 1256.3               & 59.140                & 1301.0(1.14)$\ddagger$                & 1298.5               & 0.413                           \\                                                                                 							
   
			\noalign{\smallskip}
			r106                 & -                    & -                    & \textbf{1139.5(6.29)}         & 1124.1               & 59.260               & 1257.1(29.45)$\ddagger$                 & 1205.4                 & 0.366                           \\                                                                              							

			\noalign{\smallskip}
			r107                 & 995.3                & 5164.6               & \textbf{1034.4(8.05)}         & 1025.5               & 62.400               & 1151.1(29.71)$\ddagger$                & 1085.9               & 0.380                            \\                                                                              							

			\noalign{\smallskip}
			r108                 & -                    & -                    & \textbf{ 946.5(5.93)}           & 942.0                & 61.810               & 1012.9(17.46)$\ddagger$                & 975.5                & 0.325                            \\                                                                             							

			\noalign{\smallskip}
			r109                 & 1045.7               & 13239.2              & \textbf{1086.0(12.05)}          & 1072.9               & 51.530                 & 1194.0(33.40)$\ddagger$                  & 1124.5               & 0.408                             \\                                                                           							

			\noalign{\smallskip}
			r110                 & -                    & -                    & \textbf{1042.5(5.24)}          & 1034.9               & 52.480                & 1120.8(10.02)$\ddagger$                  & 1083.3               & 0.370                             \\                                                                             							

			\noalign{\smallskip}
			r111                 & -                    & -                    & \textbf{1010.1(11.83) }        & 1004.8               & 55.735               & 1132.1(44.24)$\ddagger$                 & 1053.0               & 0.358                             \\                                                                            							

			\noalign{\smallskip}
			r112                 & 937.1                & 28734.1              & \textbf{957.6(2.41) }          & 956.4                & 55.550               & 1028.7(15.64)$\ddagger$                  & 990.8                & 0.369                            \\                                                                            							

			\noalign{\smallskip}
			w-d-l                & 6-0-6                    & -                    & -                  & -                    & -                    & 12-0-0                                   & -                    & -                                        \\
   
			\noalign{\smallskip}
			No.best              & 6                    & -                    & 0                       & 6                    & -                    & 0                                      & 0                    & -                                       \\
			\noalign{\smallskip}
			Ave.PDR              & 0.00\%                  & -                    & 3.05\%                  & -                    & -                    & 12.26\%                                  & -                    & -                                     \\
			\noalign{\smallskip}
			Ave.Time             & -                    & 12337.883            & -                       & -                    & 57.974              & -                                     & -                    & 0.380                                  \\
			\hline  \noalign{\smallskip}
		\end{tabular}
	}
\end{table*}\\
\indent In the 2LP, the currently existing algorithm  is variable neighborhood descent (VND) \cite{tagmouti2010variable}, which was as a comparison algorithm in solving 2LP. Because the vehicle departure time does not need to be calculated in the 2LP,  MAENS-GN only needs to run the algorithm in the routing stage, i.e., MAENS. For a fairer comparison,  VND* was added as a comparison algorithm, which denotes  that VND runs the same time with MAENS-GN. During the time it takes MAENS-GN to
complete one run,  VND  is run multiple times until it reaches the runtime of MAENS-GN and then stops, and finally  the best value of  the multiple values gotten by multiple runs is as the result of VND* running on that instance. 
\begin{table*}[!htb]	
	\scriptsize
	\caption{Results on the $gdb$ benchmark test set in the 3LP in terms of costs of solutions. For each instance, the average performance of an algorithm is indicated in bold if it is the best among all comparison algorithms. ``$\ddagger$" and ``$\dagger$" indicate MAENS-GN performs significantly better than and equivalently to the corresponding algorithm based on 20 independent runs according to the Wilcoxon rank-sum test with significant level $p=0.05  $, respectively.}
	\label{tab:6}
	
	\centering	
	\noindent \makebox[\textwidth][c]
	{
		\begin{tabular}{ccccccccccccc}
			\hline  \noalign{\smallskip}
			Instances     &$|R| $      & LB         &MAENS          & \multicolumn{3}{c}{MAENS-GN}        &VND   & \multicolumn{3}{c}{VND-GN}    & \multicolumn{2}{c}{VND-GN*}  \\
				\cmidrule(lr){4-4}  \cmidrule(lr){5-7}  \cmidrule(lr){8-8} \cmidrule(lr){8-8} \cmidrule(lr){9-11} \cmidrule(lr){12-13}
			&            &                      &Ave       & Ave(std)             & Best       & Time   &Ave                 & Ave(std)                  & Best       & Time      & Ave(std)   & Best\\
			\hline  \noalign{\smallskip}
			gdb1        &22       & 316         & 2510.0        & 618.5(0.90)           & 617.4      & 43.850  &2510.0              & 618.5(0.76)$\dagger$      & 617.7      & 10.310       &\textbf{617.9(0.51)} 		&617.2 \\
			\noalign{\smallskip}
			gdb2         &26      & 339        & 2132.0         & \textbf{622.4(0.92) }          & 621.2      & 65.591    &2143.8            & 631.8(12.98)$\ddagger$     & 621.1   & 12.730        &\textbf{627.6(10.27)}$\dagger$ 		&621.0\\
			\noalign{\smallskip}
			gdb3         &22      & 275         &634.0         & \textbf{289.0(0.00)}            & 289.0        & 6.703     &698.9            & 331.7(30.21)$\ddagger$    & 311.0   & 0.098           &312.2(5.45)$\ddagger$ 		&311.0\\
			\noalign{\smallskip}
			gdb4         &19     & 287       &4641.0           & \textbf{745.5(0.88)}           & 744.4       & 49.972     &4641.0           & \textbf{745.8(0.97)}$\dagger$       & 744.1   & 11.488          &\textbf{745.1(0.54)}$\dagger$ 		&744.3\\
			\noalign{\smallskip}
			gdb5         &26      & 377      &426.6            & \textbf{389.6(9.22) }          & 377.5        & 9.601   &452.8             & 416.5(9.95)$\ddagger$      & 396.5 & 0.125            &403.8(7.06)$\ddagger$ 		&391.0     \\
			\noalign{\smallskip}
			gdb6         &22      & 298      &2526.0            & \textbf{572.9(14.64)}          & 562.3        & 37.047   &2596.9            & 593.5(25.46)$\ddagger$     & 539.3    & 8.065          &\textbf{579.3(29.16)}$\dagger$ 		&538.7   \\
			\noalign{\smallskip}
			gdb7         &22      & 325      &507.2            & \textbf{468.0(3.82)}           & 459.4       & 28.250     &547.7            & 497.6(17.55)$\ddagger$     & 468.9    & 6.823          &486.7(12.32)$\ddagger$ 		&467.4   \\
			\noalign{\smallskip}
			gdb8         &46      & 348     &525.7             & \textbf{391.7(11.09)}           & 368.8      & 39.971  &550.5             & 439.5(6.43)$\ddagger$     & 411.5    & 0.655         &\textbf{397.5(13.44)}$\dagger$ 		&383.0 \\
			\noalign{\smallskip}
			gdb9        &51       & 303     &1056.5             & \textbf{319.8(10.18)}           & 309.0        & 46.583  &1145.0             & 413.2(31.61)$\ddagger$    & 348.0     & 0.899         &353.0(10.63)$\ddagger$ 		&334.0  \\
			\noalign{\smallskip}
			gdb10       &25       & 275    &340.2               & \textbf{286.6(11.19)}          & 275.0         & 8.081   &370.2             & 321.2(3.97)$\ddagger$      & 318.0     & 0.117        &313.9(7.34)$\ddagger$ 		&294.0  \\
			\noalign{\smallskip}
			gdb11       &45       & 395    &759.9              & \textbf{417.5(5.36)}          & 411.0          & 39.611  &856.6             & 468.3(17.18)$\ddagger$      & 439.0    & 0.625         &436.2(10.30)$\ddagger$ 		&417.0 \\
			\noalign{\smallskip}
			gdb12       &23       & 458     &925.4             & \textbf{479.4(25.84)}         & 458.0         & 9.007   &1020.8             & 577.6(26.76)$\ddagger$      & 540.0      & 0.088       &533.8(12.06)$\ddagger$ 		&520.0 \\
			\noalign{\smallskip}
			gdb13       &28       & 536    &705.2              & \textbf{615.1(3.39)}           & 609.8       & 31.838  &740.9              & 643.8(13.52)$\ddagger$     & 621.1    & 8.240      &641.6(14.78)$\ddagger$ 		&620.6\\
			\noalign{\smallskip}
			gdb14      &21        & 100   &342.2               & 167.5(6.13)           & 154.8       & 7.418 &353.1               & 169.8(10.18)$\dagger$      & 141.8   & 0.081      &\textbf{155.2(7.85)} 		&140.4 \\
			\noalign{\smallskip}
			gdb15      &21        & 58     &2826.0              & 79.5(1.52)             & 76.8        & 43.215 &2827.7              & 79.8(1.24)$\dagger$       & 77.2    & 5.389      &\textbf{78.2(1.06)} 		&76.2     \\
			\noalign{\smallskip}
			gdb16      &28        & 127    &152.2              & \textbf{137.1(2.61)}            & 133.1       & 11.062 &155.9              & \textbf{138.2(0.00)}$\dagger$      & 138.2     & 0.159     &\textbf{138.2(0.00)}$\dagger$ 		&138.2    \\
			\noalign{\smallskip}
			gdb17      &28        & 91     &3897.8              & 113.3(1.50)            & 109.9       & 50.624  &3900.6             & 112.7(1.26)$\dagger$      & 111.3     & 9.821     &\textbf{111.7(0.28)} 		&111.3   \\
			\noalign{\smallskip}
			gdb18      &36        & 164     &618.6             & \textbf{167.6(3.34)}           & 164.0           & 16.700  &641.0               & 192.9(12.65)$\ddagger$     & 172.0     & 0.327      &175.8(5.03)$\ddagger$ 		&166.0  \\
			\noalign{\smallskip}
			gdb19      &11        & 55       &134.0            & \textbf{60.0(0.00)}             & 60.0           & 2.045  &145.4              & 68.6(10.29)$\ddagger$     & 55.0      & 0.022     &\textbf{61.5(7.67)}$\dagger$ 		&55.0  \\
			\noalign{\smallskip}
			gdb20     &22         & 121    &269.2              & \textbf{129.4(3.92)}            & 124.0          & 7.272  &309.8               & 146.5(20.27)$\ddagger$    & 127.0    & 0.080        &\textbf{129.3(2.24)}$\dagger$ 		&126.0  \\
			\noalign{\smallskip}
			gdb21      &33        & 156     &528.1             & 258.2(11.32)            & 239.5       & 18.885 &536.9              & 262.4(3.74)$\ddagger$    & 249.6    & 0.330      &\textbf{250.1(8.09)} 		&236.9  \\
			\noalign{\smallskip}
			gdb22     &44         & 200     &593.6             & 317.7(9.57)             & 298.5       & 39.338 &599.8              & 308.2(0.14)            & 308.1   & 0.750       &\textbf{307.9(1.04)} 		&303.4 \\
			\noalign{\smallskip}
			gdb23    &55          & 233    &839.3              & 346.8(8.51)             & 335.6       & 71.849  &813.2             & 337.4(0.00)           & 337.4     & 1.690      &\textbf{337.4(0.00)} 		&337.4    \\
			\noalign{\smallskip}
			w-d-l     &-           & -     & -               & -                 & -        & -            & -                    & 15-6-2                  & -      & -             & 9-7-7            &-\\
			\noalign{\smallskip}
			No.best   &-           & -      & -              & 1                     & 18       & -            & -                    & 0                        & 5       & -            & 0            &8\\
			\noalign{\smallskip}
			Ave.PDR  &-            &-    & 638.86\%               & 37.38\%               & -       & -       & 649.94\%                & 45.22\%                   & -        & -        &39.27\%           & -\\
			\noalign{\smallskip}
			Ave.Time  &-           & -     & -               & -                      & -      & 29.761     & -          & -                        & -       & 3.431        &-               &-\\
			\hline  \noalign{\smallskip}
		\end{tabular}
	}
\end{table*}\\
\indent In the 3LP, the currently existing method is column generation (CG) \cite{tagmouti2007arc}, which serves as one comparison algorithm of MAENS-GN. The data that CG used in the comparative experiments were all provided by \cite{tagmouti2007arc}. Because  only  Solomom's test set is used in \cite{tagmouti2007arc}, CG only participated in the comparison of  Solomon's test set. Because the data \cite{tagmouti2007arc} of CG on Solomon’s test set was obtained assuming that the service time remains unchanged. Therefore, when solving Solomon’s test set, for fairness, the evaluation functions of MAENS-GN  need to be simply modified.
But CG,  being an exact algorithm,  is limited to handling small-scale instances. To evaluate  the performance of MAENS-GN on larger-scale instances,  we adapted the algorithm VND \cite{tagmouti2010variable}, which was initially designed for the 2LP. Since the primary difference between dealing with 2LP and 3LP is whether or not to optimize for vehicle departure time, the techniques  for optimizing vehicle departure time (GSS and NCS) were incorporated  into VND. VND-GN is formed in this way, which was as another comparison algorithm of MAENS-GN in the 3LP.  For a fairer comparison, VND-GN* still acted as a comparison algorithm for MAENS-GN in the 3LP, just like VND* acting as a comparison algorithm for MAENS-GN in the 2LP. VND-GN* denotes that VND-GN runs the same time as MAENS-GN.  What's more, to verify the effects of GSS and NCS, MAENS and VND were also employed  as comparison algorithms on the $ gdb $ and $ egl $ test sets in the 3LP. MAENS here denotes the method in the first stage of MAENS-GN. \\
\indent In MAENS-GN, the population size ($psize$)  was set to 10. The number of generations ($ Gm $)  was 50, and the probability of performing a local search ($Pls$) was 0.1. All other parameters in MAENS were consistent with those in \cite{tang2009memetic}. In addition, the code for VND was reproduced from the description in \cite{tagmouti2010variable}. There are two initialization methods in VND, i.e., the savings heuristic and the insertion heuristic, where the savings heuristic is a deterministic method. In the process of constructing a route in the insertion heuristic, each time one  arc is required to be selected to insert into the current route. Among all candidate arcs, the next arc selected into the route is the arc that incurs the least additional cost to the route. However, in reality, there may be multiple arcs meeting this condition simultaneously, and the article \cite{tagmouti2010variable}  does not specify the tie-breaking rule.  We opted for a stochastic tie-breaking rule, where an arc  is randomly chosen from multiple arcs meeting the condition. Therefore, the insertion heuristic is a stochastic method, and VND is thus categorized as a stochastic algorithm. Because in [9], VND solves the CARPTDSC  assuming constant service time. Therefore, the evaluation function of VND (or VND-GN) was modified here to handle the instances (except Solomon's test set) of CARPTDSC under time-dependent service time.\\
\indent The metrics are those commonly used in the literature of vehicle/arc routing~\cite{liu2021memetic,SST-2021-0372,LiuTY22,TangLYY21,LiuZTY23} and are divided into two categories. One is the metric for evaluating the solution quality, including w-d-l, No.best and Ave.PDR. The other is the time metric, namely Ave.Time \cite{liu2021memetic}. w-d-l refers to the number of instances on which the proposed algorithm is superior to, not significantly different from, or inferior to the corresponding compared algorithm in the test set. The comparison is based on the average value of 20 independent repeated experiments, and the wilcoxon ranksum test \cite{wilcoxon1992individual} with the significance level $ p =0.05 $ was employed here. No.best refers to the number of instances where the algorithm reaches the best value in the test set. The best value refers to the optimal result achieved by all comparison algorithms on that instance in 20 independent runs. Note that No.best, unless otherwise specified, refers to the data in the column headed Best in the tables, even though No.best has also data in the columns headed Ave(std). Ave.PDR refers to the value of the average performance degradation ratio (PDR)  of  the algorithm  over all instances in the instance set. PDR is calculated as $(TC_1-TC_2)/TC_2 * 100\% $, where $TC_1$ is the average cost value of the algorithm on that instance, and $TC_2$ is the lower bound  or the best value of all compared algorithms on that instance. Ave.Time refers to the average time spent by the algorithm on all instances in this test set. Therefore, the smaller the values of Ave.PDR and Ave.Time, the better. \\
\indent In Tables \ref{tab:1}-\ref{tab:15}, the columns headed Ave(std), Best, and Time stand for average results (standard deviations), best results, and average computation time (in seconds)  over 20 independent runs, respectively. In addition, the columns headed \(|R|\) provide the  number of tasks in the instances.
According to the method in \cite{tagmouti2010variable}, the lower bounds of \(gdb\) and \(egl\) in the static test set are still the lower bounds of \(gdb\) and \(egl\) in the 2LP and 3LP. Therefore the columns headed LB provide the lower bounds found so far for the instances, which were collected from \cite{tang2009memetic, tagmouti2010variable, vidal2017node, beullens2003guided, brandao2008deterministic, belenguer2003cutting, ahr2004contributions}.

\subsection{Results related to RQ1}
\label{2lp}
This section is mainly to verify whether MAENS-GN significantly improves the solution quality on the test sets in the 2LP. 

From Tables \ref{tab:1}-\ref{tab:2} and  the solution quality metrics, i.e., w-d-l, No.best and Ave.PDR, it can be seen that MAENS-GN is significantly better than VND in terms
of average values, best values and proximity to the lower
bounds on the small-scale \(gdb\) and  larger-scale \(egl\) test sets in the 2LP. For example, the w-d-l value of VND on the \(gdb\) set in the 2LP is 23-0-0, which means that the average values of MAENS-GN on all 23 instances of the \(gdb\) set in the 2LP are significantly better than those of VND. However, for  Ave.Time, MAENS-GN is more time-consuming than VND. For a fairer comparison, VND* also participated in the comparison. It can be seen from Tables \ref{tab:1}-\ref{tab:2} that when VND runs for more time, i.e., the runtime of MAENS-GN, the effect becomes better. For instance, on  the \(gdb\) set in the 2LP, the Ave.PDR values of VND and VND* are  10.82\% and 8.74\%, respectively, which is an improvement of 2.08\%. However, MAENS-GN still maintains its superiority.  The small  standard deviation  values for MAENS-GN  indicate good robustness.

\begin{sidewaystable*}[]	
	\scriptsize
	\caption{Results on the $egl$ benchmark test set in the 3LP in terms of costs of solutions. For each instance, the average performance of an algorithm is indicated in bold if it is the best among all comparison algorithms. ``$\ddagger$" and ``$\dagger$" indicate MAENS-GN performs significantly better than and equivalently to the corresponding algorithm based on 20 independent runs according to the Wilcoxon rank-sum test with significant level $p=0.05  $, respectively.}
	\label{tab:7}
	
	\centering	
	\noindent \makebox[\textwidth][c]
	{
		\begin{tabular}{ccccccccccccc}
			\hline  \noalign{\smallskip}
			Instances & $|R|$  & LB    &MAENS  & \multicolumn{3}{c}{MAENS-GN}   &VND & \multicolumn{3}{c}{VND-GN}         & \multicolumn{2}{c}{VND-GN*}  \\
				\cmidrule(lr){4-4}  \cmidrule(lr){5-7}  \cmidrule(lr){8-8} \cmidrule(lr){8-8} \cmidrule(lr){9-11} \cmidrule(lr){12-13}
			&         &        &Ave      & Ave(std)            & Best    & Time      &Ave   & Ave(std)                  & Best    & Time            & Ave(std)     & Best       \\
			\hline  \noalign{\smallskip}
			egl-e1-A &51& 3548 &3669.0 & \textbf{3588.2(7.59) }     & 3569.3  & 98.445     &3909.3 & 3867.4(0.00)$\ddagger$         & 3867.4  & 0.885    & 3858.7(17.75)$\ddagger$		&3817.7     \\
			\noalign{\smallskip}
			egl-e1-B &51& 4498 &4846.4 & \textbf{4582.2(42.20)}      & 4498.0    & 80.793  &5263.0 & 4927.8(76.50)$\ddagger$      & 4694.0    & 0.928     & 4689.2(36.57)$\ddagger$		&4628.0  \\
			\noalign{\smallskip}
			egl-e1-C &51& 5595 &6343.0 & \textbf{5895.6(62.65)}     & 5737.8  & 110.685    &6523.0 & 6125.8(21.57)$\ddagger$      & 6104.7  & 14.903     &  6092.4(36.51)$\ddagger$		&5978.7  \\
			\noalign{\smallskip}
			egl-e2-A &72& 5018 &5102.5 & \textbf{5056.4(50.04)}     & 5018.0    & 200.536  &5217.3 & 5180.0(0.00)$\ddagger$            & 5180.0    & 2.729     & 5180.0(0.00)$\ddagger$		&5180.0    \\
			\noalign{\smallskip}
			egl-e2-B &72& 6317 &6536.8 & \textbf{6377.4(22.60)}      & 6350.0    & 179.683 &6763.0 & 6605.0(0.00)$\ddagger$           & 6605.0    & 2.420      & 6577.1(38.09)$\ddagger$		&6459.0   \\
			\noalign{\smallskip}
			egl-e2-C &72& 8335 &8560.1 & \textbf{8500.3(98.31)}     & 8365.5  & 164.139  &8710.5 & 8699.0(0.00)$\ddagger$           & 8699.0    & 2.388     &  8694.4(20.05)$\ddagger$		&8607.0  \\
			\noalign{\smallskip}
			egl-e3-A &87& 5898 &6238.9 & \textbf{6091.9(99.65)}    & 5924.0    & 335.173   &6641.0 & 6425.0(48.17)$\ddagger$        & 6215.0    & 4.426        & 6380.5(71.48)$\ddagger$		&6215.0 \\
			\noalign{\smallskip}
			egl-e3-B &87& 7744 &8499.4 & \textbf{8285.4(137.10)}     & 8077.6  & 306.216 &8930.0 & 8911.9(35.90)$\ddagger$      & 8820.8  & 19.488      &  8776.2(76.53)$\ddagger$		&8576.9 \\
			\noalign{\smallskip}
			egl-e3-C &87& 10244 &10527.8 & \textbf{10404.0(61.50)}   & 10345.0   & 238.743 &11254.2 & 11109.8(151.15)$\ddagger$   & 10684.0   & 3.773      &  10699.8(68.17)$\ddagger$		&10469.0  \\
			\noalign{\smallskip}
			egl-e4-A &98& 6408 &7263.7 & \textbf{7036.9(112.66)}    & 6839.3  & 392.461  &7698.9 & 7698.9(22.23)$\ddagger$     & 7602.0    & 29.410        & 7212.6(11.57)$\ddagger$		&7163.3 \\
			\noalign{\smallskip}
			egl-e4-B &98& 8935 &9240.3 & \textbf{9190.4(57.97)}     & 9109.1  & 336.769  &9838.5 & 9805.5(109.47)$\ddagger$     & 9497.4  & 6.320         & 9545.8(62.31)$\ddagger$		&9416.8 \\
			\noalign{\smallskip}
			egl-e4-C &98& 11512 &12929.1 & \textbf{12362.4(236.91)}  & 12088.7 & 301.267 &13007.0 & \textbf{12376.3(162.53)}$\dagger$   & 12201.7 & 56.897      & 12532.7(235.41)$\ddagger$		&12038.1\\
			\noalign{\smallskip}
			egl-s1-A &75& 5018  &5193.9 & \textbf{5120.0(51.02)}       & 5058.0  & 227.150  &5603.5 & 5530.0(0.00)$\ddagger$           & 5530.0    & 3.049    & 5398.1(70.41)$\ddagger$		&5272.5    \\
			\noalign{\smallskip}
			egl-s1-B &75& 6388  &7853.6 & \textbf{7080.4(247.18)}   & 6845.0    & 251.313  &8158.6 & 7734.0(295.17)$\ddagger$      & 6789.7  & 50.189       & 7504.6(206.84)$\ddagger$		&6915.9\\
			\noalign{\smallskip}
			egl-s1-C &75& 8518  &12900.2 & 10635.6(292.52) & 10015.3 & 255.124  &12582.0 & \textbf{10318.0(475.28)}          & 10009.0   & 72.490    & 10616.9(447.63)$\dagger$		&9632.6 \\
			\noalign{\smallskip}
			egl-s2-A &147& 9825  &10778.6 & \textbf{10668.1(89.07)}  & 10454.1 & 1154.471 &11633.2 & 11621.1(92.39)$\ddagger$    & 11243.2 & 75.050        & 11265.9(103.77)$\ddagger$		&11084.7 \\
			\noalign{\smallskip}
			egl-s2-B &147& 13017 &13702.2 & \textbf{13562.7(65.58)}  & 13464.0   & 800.752  &14256.0 & 14071.0(0.00)$\ddagger$          & 14071.0   & 16.248     &  14069.1(8.28)$\ddagger$		&14033.0 \\
			\noalign{\smallskip}
			egl-s2-C &147& 16425 &19320.3 & \textbf{18521.3(402.58)} & 17700.4 & 873.643  &20029.2 & 19605.6(621.87)$\ddagger$    & 18350.2 & 117.071       & 18917.6(248.39)$\ddagger$		&1850.02\\
			\noalign{\smallskip}
			egl-s3-A &159& 10165 &10913.8 & \textbf{10901.3(68.05)}  & 10756.0   & 1456.500   &11600.8 & 11593.9(14.27)$\ddagger$     & 11556.0   & 84.768       & 11549.6(33.82)$\ddagger$		&11418.0\\
			\noalign{\smallskip}
			egl-s3-B &159& 13648 &14383.7 & \textbf{14140.5(57.27)}  & 14026.5 & 1052.030  &14834.7 & 14586.1(0.00)$\ddagger$         & 14586.1 & 25.311     & 14586.1(0.00)$\ddagger$		&14586.1   \\
			\noalign{\smallskip}
			egl-s3-C &159& 17188 &18118.1 & \textbf{18112.5(137.87)} & 17847.0   & 976.738  &19221.0 & 19220.7(146.84)$\ddagger$   & 18816.0   & 113.329      & 19019.5(132.41)$\ddagger$		&18778.0\\
			\noalign{\smallskip}
			egl-s4-A &190& 12153 &13170.0 & \textbf{12827.4(141.70)}   & 12640.0   & 1534.158 &13763.0 & 13355.0(0.00)$\ddagger$           & 13355.0   & 31.805       & 13355.0(0.00)$\ddagger$		&13355.0 \\
			\noalign{\smallskip}
			egl-s4-B &190& 16113 &17122.4 & \textbf{16831.0(122.39)} & 16643.0   & 1550.589 &17880.0 & 17647.0(0.00)$\ddagger$           & 17647.0   & 36.325       & 17531.8(131.39)$\ddagger$		&17176.0 \\
			\noalign{\smallskip}
			egl-s4-C &190& 20430 &21423.3 & 21281.8(88.73) & 21083.0   & 1380.502 &21080.8          & 20920.4(0.00)        & 20920.4  & 42.440         & \textbf{20915.0(23.69)}		&20811.7\\
			\noalign{\smallskip}
			w-d-l    &-  & -   &-       & -    & -      & -       &-         & 21-1-2                        & -       & -          & 22-1-1      &-  \\
			\noalign{\smallskip}
			No.best  &-  & -   &-       & 0        & 21      & -       &-        & 1                           & 3       & -             & 0          &3 \\
			\noalign{\smallskip}
			Ave.PDR  &-  & -   &9.50\%   & 5.74\%   & -       & -      &14.14\%           & 10.95\%                     & -       & -      &9.48\%      & - \\
			\noalign{\smallskip}
			Ave.Time &-  & -   &-       & -        & -       & 594.078    &-         & -                          & -       & 33.86         & -           &-\\
			\hline  \noalign{\smallskip}
		\end{tabular}
	}
\end{sidewaystable*}
The experiments reveal  that MAENS-GN  significantly outperforms  VND and VND* on both small-scale $ gdb $ test set  and larger-scale  $ egl $ test set for w-d-l, No.best, and Ave.PDR, suggesting that MAENS-GN significantly improves the solution quality  in the 2LP. The main reason why MAENS-GN can outperforms VND and VND* may be that the search performance of MAENS is better than that of VND. The inclusion of various search operators, such as Merge-split, helps to improve the performance of MAENS-GN. In addition, the novel initialization method makes the effect better. However, it is noted that MAENS-GN is more time-consuming than VND, primarily due to the higher time complexity of operators like Merge-split.

\subsection{Results related to RQ2}
\label{3lp}
This section mainly verifies the performance of MAENS-GN on both small-scale and larger-scale instances in the  3LP.  Additionally, the contribution of GSS and NCS to the performance of MAENS-GN is examined.

Because the comparison algorithms were run on different computers, the runtime needs to be normalized. So about the time values of CG, they were converted according to the types of computers. The computer processor used for CG in \cite{tagmouti2007arc} was the 400 MHz UltraSparc II, while the computer used for MAENS-GN and VND-GN was Intel(R) Xeon(R) Platinum 9242 CPU @ 2.30GHz. To make a fair comparison, all the time values of CG in \cite{tagmouti2007arc}    were divided by 5.75 (i.e., 2.3GHz/400MHz), then added to Tables \ref{tab:3}-\ref{tab:5}. In Table \ref{tab:5}, CG could not get results on many instances due to memory or time constraints, denoted by ``-". These instances are not involved in the calculation, such as Ave.PDR and Ave.Time.  Although CG, being an exact algorithm, outperforms MAENS-GN in terms of the solution quality (since CG can find optimal solutions), MAENS-GN significantly reduces runtime. For example, the average runtimes of CG and MAENS-GN on \(Solomon\)-25 were 626.115 and 13.875, respectively, which results in a reduction of approximately 45 times in runtime. The reason why MAENS-GN  is faster than CG is mainly because MAENS-GN is a stochastic approximate algorithm, which only needs to get some good approximate and feasible solutions. While CG is an exact algorithm, which needs to find the optimal solution, the search accuracy and time complexity are generally higher. Therefore, the search speed of MAENS-GN  is generally faster.

\begin{table*}[!htb]	
	\scriptsize
	\caption{Results on the $jd$ benchmark test set in the 3LP in terms of costs of solutions. For each instance, the average performance of an algorithm is indicated in bold if it is the best among all comparison algorithms. ``$\ddagger$" and ``$\dagger$" indicate MAENS-GN performs significantly better than and equivalently to the corresponding algorithm based on 20 independent runs according to the Wilcoxon rank-sum test with significant level $p=0.05  $, respectively.}	
	\label{tab:8}
	\centering	
	\noindent \makebox[\textwidth][c]
	{
		\begin{tabular}{cccccccc}
			\hline  \noalign{\smallskip}
			Instances    & $|R|$        & \multicolumn{3}{c}{MAENS-GN}                                                       & \multicolumn{3}{c}{VND-GN}                                                                         \\
				\cmidrule(lr){3-5}  \cmidrule(lr){6-8}  
			            &        & Ave(std)              & Best                 & Time                 & Ave(std)                   & Best                 & Time                                        \\
			\hline  \noalign{\smallskip}
			jd-50-1     &50       & \textbf{3268.2(11.80)}        & 3252.8               & 79.554               & 3386.7(18.44)$\ddagger$             & 3334.1               & 1.872                                \\
			\noalign{\smallskip}
			jd-50-2     &50       & \textbf{3593.8(32.80)}          & 3546.0                 & 87.707               & 3690.6(40.37)$\ddagger$             & 3612.0                 & 10.185                             \\
			\noalign{\smallskip}
			jd-50-3     &50       & \textbf{3363.5(8.92)}          & 3342.0                 & 51.834               & 3484.1(37.69)$\ddagger$             & 3423.0                 & 1.677                              \\
			\noalign{\smallskip}
			jd-50-4    &50        & \textbf{3580.2(28.08)}        & 3514.1               & 50.784               & 3710.0(38.18)$\ddagger$             & 3672.0                 & 9.003                               \\
			\noalign{\smallskip}
			jd-100-1    &100       & \textbf{6234.8(18.29)}        & 6211.0                 & 569.377              & 6411.4(47.80)$\ddagger$              & 6348.0                 & 18.275                            \\
			\noalign{\smallskip}
			jd-100-2    &100       & \textbf{6186.2(25.71)}          & 6117.0                 & 555.759              & 6383.2(55.38)$\ddagger$            & 6308.0                 & 17.363                            \\
			\noalign{\smallskip}
			jd-100-3    &100       & \textbf{5876.3(33.83)}          & 5822.2               & 561.251              & 6194.8(48.66)$\ddagger$             & 6088.2               & 13.642                             \\
			\noalign{\smallskip}
			jd-100-4    &100       & \textbf{6280.4(33.85)}           & 6213.0                 & 552.321              & 6467.9(67.47)$\ddagger$           & 6317.0                 & 16.815                              \\
			\noalign{\smallskip}
			jd-150-1   &150        & \textbf{8992.4(36.63)}          & 8916.9               & 1584.547             & 9316.4(74.51)$\ddagger$             & 9186.0                 & 97.842                              \\
			\noalign{\smallskip}
			jd-150-2   &150        & \textbf{9298.0(56.46)}            & 9188.0                 & 1598.263             & 9712.7(54.72)$\ddagger$            & 9605.0                 & 106.934                              \\
			\noalign{\smallskip}
			jd-150-3   &150        & \textbf{9277.6(63.22)}          & 9131.0                 & 1544.837             & 9644.1(84.76)$\ddagger$              & 9453.6               & 95.801                              \\
			\noalign{\smallskip}
			jd-150-4   &150        & \textbf{9188.0(40.36)}            & 9080.0                 & 1479.287             & 9380.8(88.57)$\ddagger$            & 9205.0                 & 62.927                             \\
			\noalign{\smallskip}
			jd-200-1   &200        & \textbf{10865.5(52.49)}           & 10755.4              & 3260.699             & 11442.5(116.23)$\ddagger$        & 11074.9              & 114.503                           \\
			\noalign{\smallskip}
			jd-200-2   &200        & \textbf{10866.2(45.73)}         & 10770.1              & 3311.978             & 11351.7(72.46)$\ddagger$          & 11190.3              & 122.243                            \\
			\noalign{\smallskip}
			jd-200-3   &200        & \textbf{11727.0(35.86)}            & 11658.0                & 3018.696             & 12117.7(124.04)$\ddagger$        & 11870.6              & 225.431                          \\
			\noalign{\smallskip}
			jd-200-4   &200       & \textbf{11529.8(52.23)}        & 11440.0                & 3052.575             & 11767.8(107.15)$\ddagger$        & 11593.0                & 135.499                        \\
			\noalign{\smallskip}
			w-d-l      &-          & -                   & -                    & -                      & 16-0-0                  & -                    & -                                           \\
			\noalign{\smallskip}
			No.best    &-          & 0                        & 16                   & -                    & 0                       & 0                    & -                                        \\
			\noalign{\smallskip}
			Ave.PDR    &-         & 0.00\%                   & -                    & -                    & 3.58\%                   & -                    & -                                        \\
			\noalign{\smallskip}
			Ave.Time   &-          & -                         & -                    & 1334.967             & -                       & -                    & 65.626                                   \\
			\hline  \noalign{\smallskip}
		\end{tabular}
	}
\end{table*}

As can be seen from Tables \ref{tab:3}-\ref{tab:8}, in terms of solution quality metrics w-d-l, No.best, and Ave.PDR, MAENS-GN is significantly better than VND-GN in the \(Solomon\)-25, \(Solomon\)-35, \(Solomon\)-40, \(gdb\), \(egl\) and, real-world \(jd\) test sets in the 3LP. For example, on the \(jd\)  set, the w-d-l value of VND-GN is 16-0-0, indicating that the average values of MAENS-GN on all 16 instances in the \(jd\)  set are significantly better than those of VND-GN. However on Ave.Time, MAENS-GN is  worse than VND-GN. For a fairer comparison, VND-GN* is compared on the small-scale \(gdb\) test set  and the larger-scale \(egl\) test set  in the 3LP in Tables \ref{tab:6}-\ref{tab:7}. It can be seen that as the runtime increases, the performance  of VND-GN improves. For example, on the \(gdb\) set in the 3LP in Table \ref{tab:6}, the Ave.PDR values of VND-GN and VND-GN* are 45.22\% and 39.27\%, respectively, which indicates an improvement of about 6\%. Despite this improvement, MAENS-GN remains  better than VND-GN*. \\
\indent In the comparison between MAENS and MAENS-GN, as well as  between VND and VND-GN, presented in Tables \ref{tab:6}-\ref{tab:7}, it is observed  that GSS and NCS have greatly improved the performance  of both MAENS and VND. For example, it can be seen from Table \ref{tab:6} that the Ave.PDR values of MAENS and MAENS-GN are 638.86\% and 37.38\%, respectively, which indicates that GSS and NCS improved MAENS by about 17 times. It can also be seen that the difference between MAENS-GN and VND-GN  is almost equal to the difference between MAENS and VND on Ave.PDR. For example, on the \(egl\) set in the 3LP in Table \ref{tab:7}, the Ave.PDR values of VND-GN and MAENS-GN are 10.95\% and 5.74\%, respectively, and the difference is 5.21\%. The  Ave.PDR values of VND and MAENS  are 14.14\% and 9.50\%, respectively, and the difference is 4.64\%. This shows that the reason why MAENS-GN is better than VND-GN is largely due to the fact that MAENS is better than VND.\\
\begin{table*}[!htb]	
	\scriptsize
	\caption{Comparison between MAENS* and MAENS on the $egl$ benchmark test set in the 3LP in terms of costs of solutions.}	
	\label{tab:15}
	\centering	
	\noindent \makebox[\textwidth][c]
	{
		\begin{tabular}{cccccccccc}
			\hline  \noalign{\smallskip}
			Instances& $|R|$ & LB    & MAENS*   & MAENS              & \multicolumn{1}{|c}{Instances}      & $|R|$ & LB    & MAENS*    & MAENS   \\
			\hline  \noalign{\smallskip}
			egl-e1-A &51 & 3548  & 3550.6  & \textbf{3549.2}  & \multicolumn{1}{|c}{egl-s1-A} &75 & 5018  & \textbf{5145.9}   & 5178.2   \\
			\noalign{\smallskip}
			egl-e1-B &51 & 4498  & 4576.9  & \textbf{4575.7}  & \multicolumn{1}{|c}{egl-s1-B} &75 & 6388  & \textbf{6576.6}   & \textbf{6576.6}   \\
			\noalign{\smallskip}
			egl-e1-C &51 & 5595  & \textbf{5677.6}  & 5679.8  & \multicolumn{1}{|c}{egl-s1-C} &75 & 8518  & \textbf{8783.8}   & 8803.6   \\
			\noalign{\smallskip}
			egl-e2-A &72 & 5018  & \textbf{5102.4}  & 5116.9  & \multicolumn{1}{|c}{egl-s2-A} &147 & 9825  & \textbf{10490.1}  & 10514.6  \\
			\noalign{\smallskip}
			egl-e2-B &72 & 6317  & \textbf{6424.9}  & 6442.9  & \multicolumn{1}{|c}{egl-s2-B} &147 & 13017 & \textbf{13887.2}  & 13903.4  \\
			\noalign{\smallskip}
			egl-e2-C &72 & 8335  & 8599.0    & \textbf{8597.1}  & \multicolumn{1}{|c}{egl-s2-C} &147 & 16425 & \textbf{17392.6}  & 17399.6  \\
			\noalign{\smallskip}
			egl-e3-A &87 & 5898  & 6099.4  & \textbf{6047.4}  & \multicolumn{1}{|c}{egl-s3-A} &159 & 10165 & \textbf{10733.0}    & 10800.9  \\
			\noalign{\smallskip}
			egl-e3-B &87 & 7744  & 8040.6  & \textbf{8032.0}    & \multicolumn{1}{|c}{egl-s3-B} &159 & 13648 & \textbf{14418.0}    & 14490.4  \\
			\noalign{\smallskip}
			egl-e3-C &87 & 10244 & \textbf{10460.6} & 10462.4 & \multicolumn{1}{|c}{egl-s3-C} &159 & 17188 & \textbf{18269.3}  & 18279.0    \\
			\noalign{\smallskip}
			egl-e4-A &98 & 6408  & \textbf{6710.6}  & 6729.0    & \multicolumn{1}{|c}{egl-s4-A} &190 & 12153 & 13197.1  & \textbf{13193.5}  \\
			\noalign{\smallskip}
			egl-e4-B &98 & 8935  & \textbf{9329.9}  & 9337.2  & \multicolumn{1}{|c}{egl-s4-B} &190 & 16113 & \textbf{17260.1}  & 17302.6  \\
			\noalign{\smallskip}
			egl-e4-C &98 & 11512 & \textbf{11992.6} & 12004.3 & \multicolumn{1}{|c}{egl-s4-C} &190 & 20430 & \textbf{21796.0}    & 21880.0    \\
			\noalign{\smallskip}
			Ave.CD  &-   & -     & -       & -       &-   &-   & -     & 482.28 & 498.18 \\
			\noalign{\smallskip}
			Ave.PDR &-   & -     & -       & -       &-   &-  & -     & 4.18\%   & 4.32\%  \\
			\hline  \noalign{\smallskip}
		\end{tabular}
	}
\end{table*}
\indent In summary, MAENS-GN demonstrates excellent results across  the small-scale, larger-scale, or real-world test sets in the 3LP. Moreover, GSS and NCS have greatly enhanced  the performance of MAENS-GN.

\subsection{Effectiveness of the novel initialization strategy}
\label{init}
This section aims to investigate  a question: Is the novel initialization strategy in MAENS superior to  the original initialization strategy? Here, two algorithms are compared, namely MAENS with the original initialization strategy (here represented as MAENS) and MAENS with the novel initialization strategy (here represented as MAENS*). The larger-scale $ egl $ test set  in the 2LP was used.  Ave.CD and Ave.PDR were used as the metrics on the test set. Ave.CD represents the average CD on all instances. CD, representing  the cost difference from the lower bound (i.e., LB) on an instance, is expressed as $CD=(cost-LB).$  \\
\indent In Table \ref{tab:15}, the Ave.CD values of MAENS* and MAENS are 482.28 and 498.18, respectively, with a difference of 15.9. This observation indicates that the novel initialization strategy enhances  the cost value of the algorithm by 15.9 on each instance on average. Turning to Ave.PDR, the values of MAENS* and MAENS are 4.18\% and 4.32\%, respectively, with a difference of 0.14\%. This  shows that the novel initialization strategy improves the algorithm  by 0.14\% on Ave.PDR. \\
\indent Through the aforementioned  experiments, it is evident that, on the larger-scale $ egl $ test set in the 2LP, the novel initialization strategy has visible  improvements to the proposed algorithm in terms of the average costs of the solutions. However, experiments, conducted on the small-scale  $ gdb $ test set in the 2LP, did not exhibit substantial improvement with  the novel initialization strategy, and the results are not presented  here. The main reason behind this observation may be that,  as the number of the tasks increases, the quality of the initialized solutions  deteriorates when using the original initialization strategy.
MAENS with the original initialization strategy already demonstrates  satisfactory  results on the small-scale \(gdb\) test set in the 2LP, making  the impact  of the operation for the time-dependent characteristics less apparent. However, on the larger-scale $ egl $ test set, the effectiveness  of the original initialization strategy diminishes. At this time, the incorporation  of targeted operations for time-dependent characteristics proves to enhance the effect of the initialization strategy. It is worth mentioning why experiments were not conducted on the test sets in the 3LP, such as the \(gdb\) and \(egl\) sets in the 3LP. One reason is due to space limitations.
Of course, the main reason is that the original/novel initialization strategy is embedded in the algorithm in the first stage of MAENS-GN, i.e., MAENS, and assessing the effect of MAENS provides a more direct perspective. Additionally, as indicated in Section \ref{3lp}, the superiority of MAENS-GN over VND-GN is primarily attributed to the superiority of MAENS over VND. This insight suggests that an improvement in the effect of MAENS generally leads to the improvement of performance of MAENS-GN on the test sets in the 3LP. Therefore, it is not really necessary to conduct experiments on the test sets in the 3LP either.
\section{Conclusion and future work}
\label{concf}
Experiments were  conducted on the modified standard test sets and a real-world test set, which demonstrates  excellent  results for MAENS-GN. Specifically, the proposed algorithm, i.e., MAENS-GN, has effectively  addressed two research questions, which indicates that  two existing limitations are addressed in this field.\\
\indent In the optimization of the vehicle departure time, this paper conducts a detailed analysis of the relationship between the route cost and the vehicle departure time in various scenarios. Based on the analyzed relationship, tailored  methods are employed  in different scenarios to optimize  the vehicle departure time, which thereby avoids the waste of calculation time and lack of accuracy of the obtained departure time caused by using a uniform  method. Furthermore, in the routing stage, a specialized  operation  for time-dependent characteristics is introduced to the initialization strategy to improve its effectiveness.\\
\indent Despite its significant advantages, MAENS-GN also has some shortcomings.
Firstly, MAENS-GN exhibits  high time consumption.
Secondly, although many solutions obtained by MAENS-GN  reach the lower bounds of the instances in the 2LP,  there remains a  discernible gap from the lower bounds of the instances in the 3LP.
Consequently, our next work focuses on  reducing the time consumption of the algorithm and enhancing  the solution quality in the 3LP.
One way to improve solution quality in the 3LP can be as follows. With the reduction of time consumption of MAENS-GN, the algorithms in the two stages of MAENS-GN, i.e., MAENS and GSS (or NCS),  can be iterated multiple times to  further enhance  the quality of the solutions in the 3LP.
\bibliographystyle{elsarticle-num} 
\bibliography{ref}





\end{document}